\documentclass[11pt]{article}

\usepackage{graphicx}
\usepackage[top=40pt, bottom=50pt, left=50pt, right=50pt]{geometry}
\usepackage[english]{babel}
\usepackage{amsmath,amssymb}
\usepackage{epstopdf}
\usepackage{tikz,url,subfig,listings,verbatim,hyperref}
\usepackage[labelfont=bf,justification=centering]{caption}
\usepackage{cleveref} % The cleveref package should automatically give Figure 1, (2), Table 5 etc: see \captionsenglish

    \crefname{figure}{Figure}{Figures}
    \Crefname{figure}{Figure}{Figures}
    \crefname{table}{Table}{Tables}
    \Crefname{table}{Table}{Tables}
    \crefname{section}{\S}{\S}
    \Crefname{section}{\S}{\S}
    \crefname{equation}{}{}
    \Crefname{equation}{}{}
    \crefname{remark}{Remark}{Remarks}
    \Crefname{remark}{Remark}{Remarks}

\newcommand{\MK}{{\mathcal K}}
\newcommand{\MO}{{\mathcal O}}

\DeclareMathOperator{\SUPP}{supp}

\usetikzlibrary{patterns,arrows,decorations.markings}

\usetikzlibrary{decorations.pathreplacing,calc,shapes}

\tikzset{middlearrow/.style 2 args={
        decoration={markings,
            mark= at position #1 with {\arrow{#2}} ,
        },        postaction={decorate}    }	}

\begin{document}

\title{High-frequency asymptotic compression of dense BEM matrices for general geometries without ray tracing}

\author{
Daan Huybrechs\footnote{daan.huybrechs@cs.kuleuven.be}\\
Department of Computer Science\\
KU Leuven, Belgium
\and
Peter Opsomer\footnote{peter.opsomer@cs.kuleuven.be}\\
Department of Computer Science\\
KU Leuven, Belgium
}
\maketitle

\begin{abstract}
Wave propagation and scattering problems in acoustics are often solved with boundary element methods. They lead to a discretization matrix that is typically dense and large: its size and condition number grow with increasing frequency. Yet, high frequency scattering problems are intrinsically local in nature, which is well represented by highly localized rays bouncing around. Asymptotic methods can be used to reduce the size of the linear system, even making it frequency independent, by explicitly extracting the oscillatory properties from the solution using ray tracing or analogous techniques. However, ray tracing becomes expensive or even intractable in the presence of (multiple) scattering obstacles with complicated geometries. In this paper, we start from the same discretization that constructs the fully resolved large and dense matrix, and achieve asymptotic compression by explicitly localizing the Green's function instead. This results in a large but sparse matrix, with a faster associated matrix-vector product and, as numerical experiments indicate, a much improved condition number. Though an appropriate localisation of the Green's function also depends on asymptotic information unavailable for general geometries, we can construct it adaptively in a frequency sweep from small to large frequencies in a way which automatically takes into account a general incident wave. We show that the approach is robust with respect to non-convex, multiple and even near-trapping domains, though the compression rate is clearly lower in the latter case. Furthermore, in spite of its asymptotic nature, the method is robust with respect to low-order discretizations such as piecewise constants, linears or cubics, commonly used in applications. On the other hand, we do not decrease the total number of degrees of freedom compared to a conventional classical discretization. The combination of the sparsifying modification of the Green's function with other accelerating schemes, such as the fast multipole method, appears possible in principle and is a future research topic.
%\keywords{Boundary Element Method \and Oscillatory integration \and High-frequency scattering \and Compression \and Condition number \and Smooth window functions}
%\subclass{65N38 \and 65F50 \and 45A05 \and 45M05 \and 65R20}
\end{abstract}

\section{Introduction} 

Numerical simulations in acoustics are often based on a Boundary Integral equation reformulation of the Helmholtz equation, see for example \cite{ColtonKress,BEA}. An incoming wave that is scattered by an obstacle with boundary $\Gamma$ results in a scattered field $u^{\text{s}}(\mathbf{x})$ that can be represented by the so-called \emph{single layer potential},
\begin{equation}
 u^{\text{s}}(\mathbf{x}) = \int_{\Gamma} K(\mathbf{x},\mathbf{y}) v(\mathbf{y}) {\rm d}s(\mathbf{y}). \label{Eslpot}
\end{equation}
Here, $K(\mathbf{x},\mathbf{y})$ is the Green's function of the Helmholtz equation with wavenumber $k$ and $v(\mathbf{y})$ is the unknown \emph{density function} defined on $\Gamma$. It is a physical quantity, namely the normal derivative of the total wave. Sound-soft and sound-hard obstacles give rise to Dirichlet or Neumann boundary conditions respectively, with zero pressure or zero normal velocity on $\Gamma$. Integral representations of the scattered field such as the one above, coupled with a boundary condition, give rise to an integral equation. In this paper, we present our results for the following integral equation of the first kind with a Dirichlet boundary condition:
\begin{equation}
 \int_{\Gamma} K(\mathbf{x},\mathbf{y}) v(\mathbf{y}) {\rm d}s(\mathbf{y}) = -u^{\text{inc}}(\mathbf{x}), \quad \forall \mathbf{x} \in \Gamma. \label{Eie} 
\end{equation}
Other representations of the scattered field and integral equations exist, such as the combined double and single layer potential \cite[\S 3.6]{ColtonKress} and the Brakhage-Werner integral equation \cite{BrakhageWerner}. These also depend on the boundary condition and on whether we consider the interior or exterior problem, but the principle of our method applies similarly as long as the integral operator asymptotically localizes.

The function $u^{\text{inc}}(\mathbf{x})$ in the right hand side of \eqref{Eie} represents an incoming wave, and the integral equation expresses that the total field $u^{\text{tot}}(\mathbf{x}) = u^{\text{s}}(\mathbf{x}) + u^{\text{inc}}(\mathbf{x})$ vanishes on the boundary of the domain. The discretization of $v(\mathbf{y})$ in \eqref{Eie} results in a boundary element method. 

As we shall consider high wavenumbers, this density function $v(y)$ is highly oscillatory, thus requiring many degrees of freedom. In order to resolve all the oscillations of the problem, we follow a standard rule of thumb and choose a number of degrees of freedom $N$ that is proportional to the wavenumber $k$ in 2D, and proportional to $k^2$ in 3D. A collocation or Galerkin discretization of \eqref{Eie} results in a linear system
\begin{equation}
\label{Esystem}
Ac=b,
\end{equation}
where $b$ represents the boundary condition. For both kinds of discretizations, the matrix $A$ is dense, large and increasingly ill-conditioned. The relationship between the condition number and the wavenumber is complicated, and depends on properties of the obstacle and of the integral operators, see~\cite[\S 5]{reviewLangdon} and \cite{condKress}. It should be noted that integral equations of the second kind are more popular in applications because they are typically better conditioned for increasing $N$. However, the condition number also deteriorates when increasing the wavenumber.

A fast matrix-vector product with $A$ in the high-frequency regime is achieved by the high-frequency Fast Multipole Method \cite{FMM,greengard1987fmm,widebandFMM,Ying2015FastDir}. Alternatively, so-called hybrid numerical-asymptotic methods aim to significantly reduce the size of the linear system by incorporating information about the solution from asymptotic analysis (see the review \cite{reviewLangdon} and references therein). In particular, phase-extraction methods use information about the phase $g$ of the solution $v(\mathbf{y})$ in order to discretize only the remaining non-oscillatory part $f$ in the following factorization:
\begin{equation}
 v(\mathbf{y}) = f(\mathbf{y},k)e^{ikg(\mathbf{y})}. \label{Eextraction}
\end{equation}
Phase extraction methods are simplest for convex obstacles, and require ray-tracing or similar techniques for more complicated domains or multiple scattering configurations \cite{bruno,geuzaine2005multiple,2DEcevit,3DEcevit,groth,chandlerwilde2015nonconvex}. These methods typically lead to matrices that are small, sometimes with frequency-independent size, but that are still dense. The use of specialized quadrature methods for highly oscillatory integrals in \cite{sparseDiscr} led to a matrix that is both small and sparse. It is possible to compute asymptotic information (more specifically, stationary points) by computing where derivatives or gradients of the phase vanish, see \cite{Ganesh2}. This is not unlike ray tracing.

The goal of this paper is to explore the use of asymptotic properties of the integral operator, but without performing phase extraction or ray tracing techniques. In doing so, we aim for great generality and flexibility in the shape of the domain. We use no a priori asymptotic information about the solution, nor do we numerically estimate the phase of $v(\mathbf{y})$. For complicated domains, the solution could be a sum of many functions of the form \eqref{Eextraction}, or it could have more complicated asymptotic structure, and that limits applicability of purely asymptotic methods.

The methodology of our approach is as follows. The linear system \eqref{Esystem} represents the action of an oscillatory integral operator (discretized in $A$) on an oscillatory function (discretized in $c$). Hence, the product of a row of $A$ with the column vector $c$ discretizes a highly oscillatory integral. Without ray-tracing methods and phase estimation, the phase of the oscillations of this integral is unknown. Yet, each highly oscillatory integral of this kind is essentially local in nature, and one can localize oscillatory integrals explicitly using smooth cut-off functions. Conceptually, in our setting this corresponds to modifying the Green's function, or more generally the kernel of the integral equation at hand. That is, we replace $K(\mathbf{x},\mathbf{y})$ by
\begin{equation}
 \tilde{K}(\mathbf{x},\mathbf{y}) = K(\mathbf{x},\mathbf{y}) w(\mathbf{x},\mathbf{y}),\nonumber
\end{equation}
where $w(\mathbf{x},\mathbf{y})$ is a smooth window function, or a sum of smooth window functions. This results in a modified discretization matrix $\tilde{A}$. The integral operator with kernel $\tilde{K}$ agrees asymptotically with the original operator as long as $w(\mathbf{x},\mathbf{y})$ equals $1$ in regions that contribute to the oscillatory integrals, and smoothly decays to zero outside these regions. Using basis functions with compact support in the discretization, the corresponding entries of $\tilde{A}$ are exactly equal to those of $A$. For regions where $w(\mathbf{x},\mathbf{y})$ vanishes, the entries of $\tilde{A}$ are exactly zero. In the intermediate regions, the entries of $\tilde{A}$ are perturbations of those of $A$.

The idea and implementation of this approach is fairly straightforward, and can be based on an existing BEM implementation. We have performed our computations in Matlab for 2D problems. Our implementation is publicly available on GitHub \cite{github} and a copy of the code can always be obtained from the authors. The main challenge is estimating suitable locations for the window functions. In principle, this still requires asymptotic analysis of the scattering problem, and this remains complicated for complicated geometries. However, the task of locating windows is significantly simpler than the task of determining (multiple) phases of the solution. Indeed, some contributions to the oscillatory integrals intrinsically correspond to reflections of rays in the scattering configuration.

Exploiting a simple geometric visibility criterion results in some compression for any obstacle, or indeed any combination of obstacles in a multiple scattering configuration, regardless of the geometric shape. But the compression rate can be much improved still. We obtain the best results by numerically measuring the contibuting regions of the integrals at a moderate value of $k$, and extrapolating the location of the windows to a higher value of $k$. Their location is refined with every significant increase of $k$, resulting in higher compression rates for high-frequency problems in an overall fairly efficient frequency sweeping solver. This does not quite achieve the effectiveness of a true asymptotic scheme in problems where those apply, but our scheme is applicable to a much wider range of problems, and the results always improve with increasing wavenumber $k$.

Sparsifying $A$ results in lower CPU times, especially in a frequency sweep. Another important observation is that the condition number of $\tilde{A}$ seems much lower than that of $A$ itself. The lower condition number has a positive effect on the number of iterations required in an iterative solver for \eqref{Esystem}. Furthermore, we show that the method is robust in several ways. It is robust with respect to the geometry of the obstacle: some amount of compression can still be achieved even for near-trapping domains, which are a worst case for ray-tracing methods. It is also robust with respect to the order of the BEM discretization. In spite of the fact that the asymptotic analysis of oscillatory integrals is mostly based on $C^\infty$ integrands, we show that asymptotic compression can be achieved even when $v(\mathbf{y})$ is discretized using piecewise linear or cubic functions. In this case, the oscillatory integrals that we localize are merely continuous. Exploiting this robustness, we can approximate the matrix entries of $\tilde{A}$ by a simple rescaling of the entries of $A$ that are computed in an existing code.

On the other hand, unlike ray tracing methods, our method is not frequency independent. Also, similar to asymptotic methods, many parts of our calculations depend on the nature and direction of the incident wave. However, the sparsifying modification of the Green's function is compatible with the Fast Multipole Method, since all we do is modify the Green's function, and a more efficient implementation in 3D exploiting that connection is a possible future research project.

The structure of the paper is as follows. We review asymptotic properties of oscillatory integrals in \cref{s:integrals} and one can find the asymptotic compression scheme in \cref{Smethod}. The localisation of the smooth functions using adaptive computation of correlations and adaptive recompression for increasing frequency is described in \cref{s:adaptive}. Finally, we end the paper with a series of numerical experiments on a range of scattering configurations in \cref{Snum} and with some concluding remarks in \cref{Sconc}.

\section{Oscillatory integrals}
\label{s:integrals}

\subsection{Asymptotic analysis of oscillatory integrals}
\label{Sintegrals}

Oscillatory integrals of the form
\begin{equation}
 \int_a^b f(x) e^{i k g(x)} {\rm d}x \label{eq:oscillatoryintegral}
\end{equation}
are a classical topic in asymptotic analysis \cite{BleisteinHandelsman,Wong}. Note that in our setting, the total phase of the integrand $g(x)$ is the sum of the phase of the Green's function, which is known explicitly in general, and the phase (or multiple phases) of a density function $v$ of the form \cref{Eextraction}. The main observation is that the integral can be seen as a sum of contributions, originating in a small number of critical points $p$. They include the endpoints of integration ($a$ and $b$ in the integral above), any singularities of the integrand and so-called stationary points. The latter are points $p$ where the derivative of the phase function $g(x)$ vanishes, i.e. $g'(p)=0$. Intuitively, the explanation is that near a stationary point the integrand is locally non-oscillatory, but elsewhere oscillations cancel each other out. Similarly, the cancellation effect is reduced or absent near the endpoints of the interval and near singularities.

Mathematically, the result of the analysis is an asymptotic expansion of the form
\[
 \int_a^b f(x) e^{i k g(x)} {\rm d}x \sim \sum_{m=1}^M e^{i k g(p_m)} \sum_{j=1}^\infty A_{m,j}[f,g] k^{-r_{m,j}}, \qquad k \to \infty.
\]
Here, the points $p_m$ are all the critical points of the integral, and the exponents $r_{m,j}$ form a strictly increasing sequence for each $m$. The coefficients $A_{m,j}$ depend linearly on $f$ and non-linearly on $g$. Explicit expressions can be derived, showing that $A_{m,j}$ depends on the function values of $f$, $g$ and their derivatives at the critical point $p_m$. The larger the index $j$, the higher the order of the derivative involved. Note that the expansion above assumes that $f$ and $g$ are infinitely differentiable. The exponents $r_{m,j}$ depend on the nature of the contributing point, and they are typically integer or rational numbers. We refer the reader to \cite{deano2017book} for a detailed introduction to asymptotics of oscillatory integrals, examples of the computation of such expansions and numerical methods to evaluate such integrals efficiently for large frequency $k$.

Asymptotic expansions can be obtained for similar integrals involving an oscillatory Hankel function, for example, rather than the complex exponential function. Since the Green's function of the Helmholtz equation in 2D is given in terms of the Hankel function, these expansions can in fact be used to deduce the asymptotic expansions of the solutions of certain integral equations \cite{groth2017highorder}. Yet, we are not concerned with the computation of asymptotic expansions in this paper. Instead, we solely rely on the fact that they exist and exploit the localization of oscillatory integrals by more robust means.

\subsection{Explicit localization of oscillatory integrals}

\begin{figure}[t]
\centering
\subfloat{\includegraphics[width=0.5\hsize]{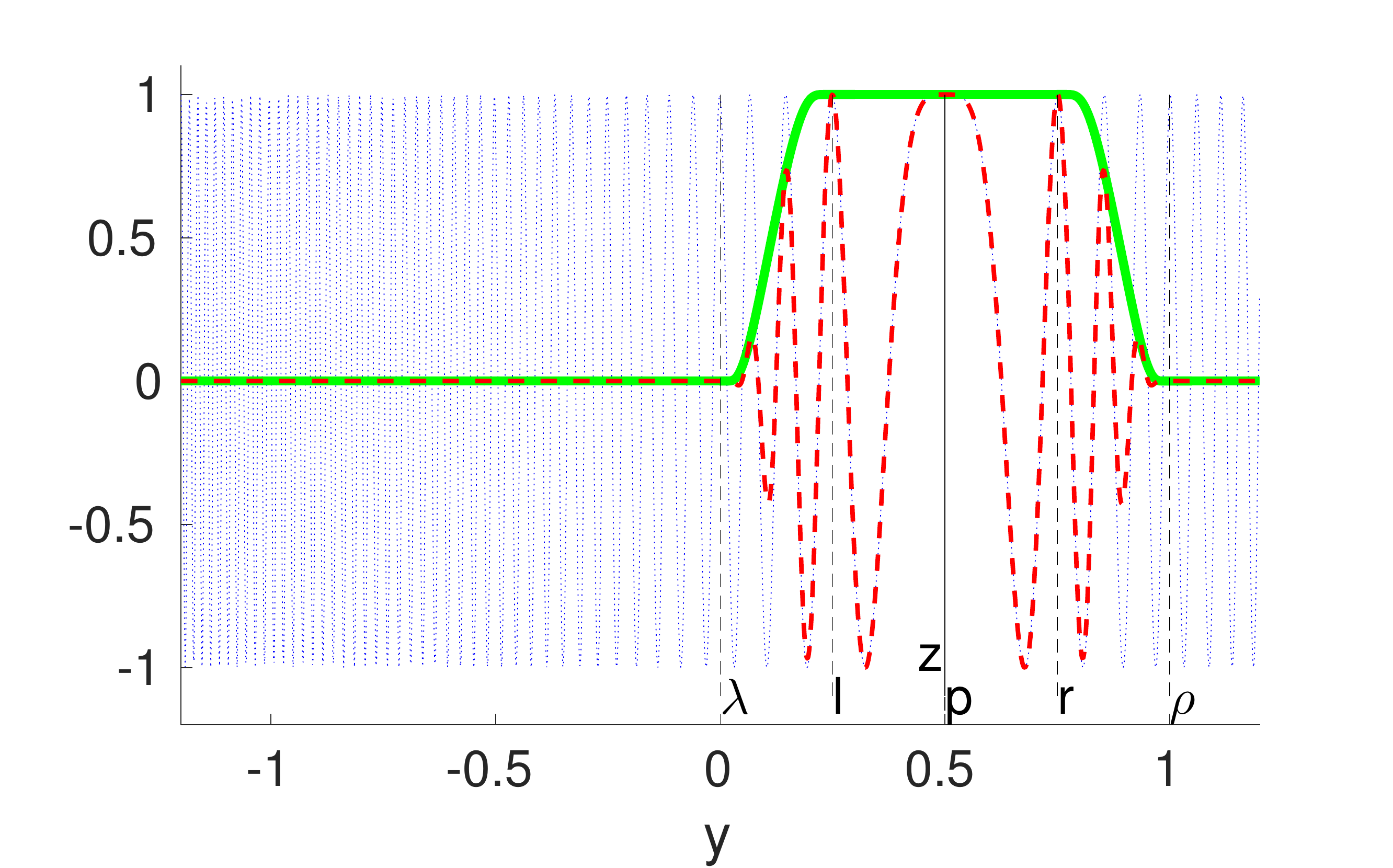} }
\subfloat{\includegraphics[width=0.5\textwidth]{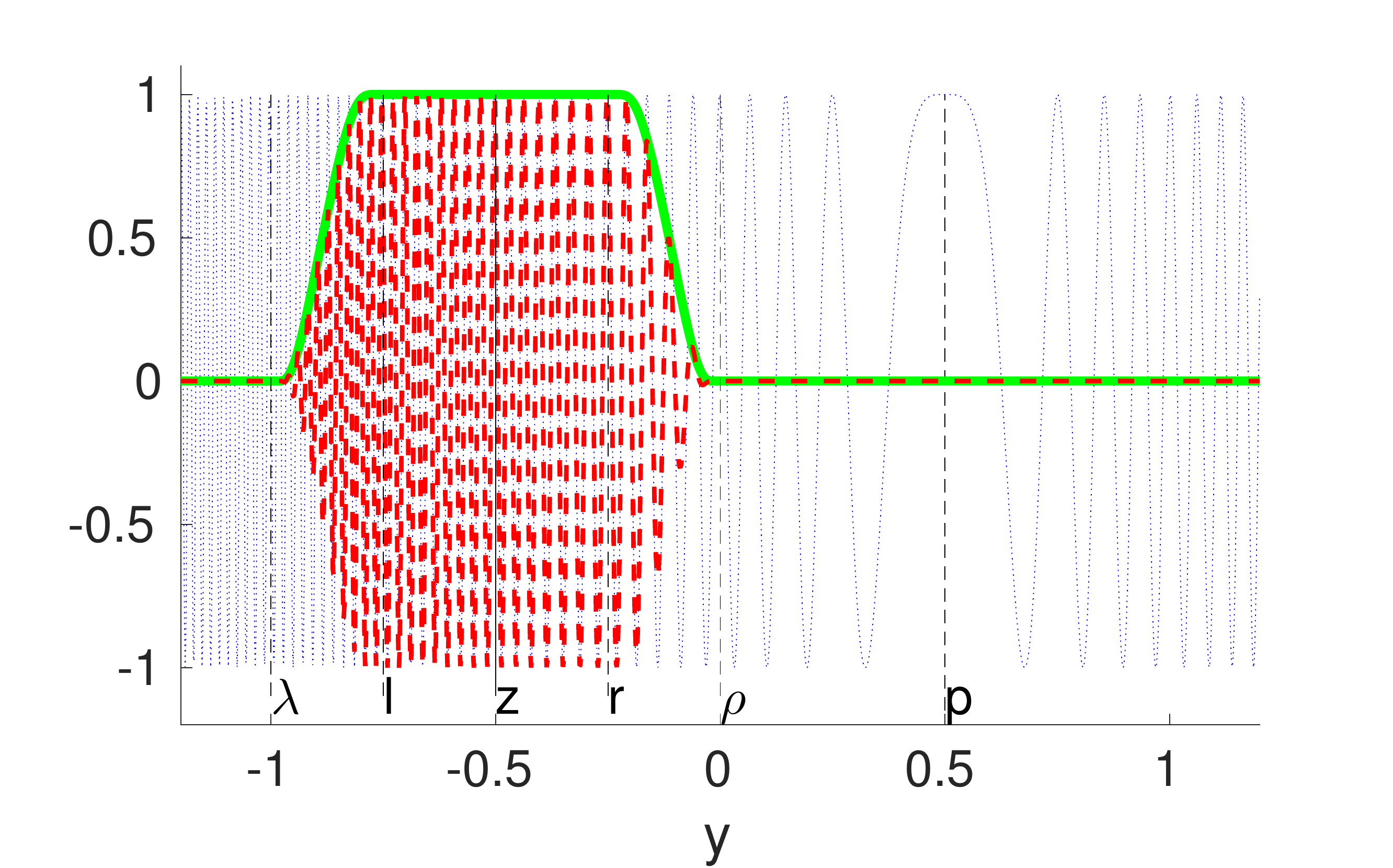} }
\caption{Illustration of the real part of an oscillatory integrand $e^{i k (y-1/2)^2}$ (dotted blue) with $k=100$ and with a stationary point at $y=1/2$. The stationary point contribution can be computed using a smooth window function (solid green) that covers it. The integral of the windowed function (dashed red) in the left panel asymptotically agrees with the integral of the oscillatory function. In contrast, the integral of the windowed function in the right panel is superalgebraically small in $k$.}
\label{FoscIntIll}
\end{figure}

It is well-known, and common in applications, to single out the contributions of oscillatory integrals using smooth cut-off functions (see, for example, \cite{bruno}). This is illustrated in \cref{FoscIntIll}, which shows the real part of the oscillatory function $e^{ik(y-1/2)^2}$ with $k=100$. This is a simple case of the general oscillatory integral \eqref{eq:oscillatoryintegral}, with a stationary point at $p=1/2$. The oscillations rapidly cancel out for high $k$, so that we can neglect all but a small region around the stationary point. We can evaluate this contribution by multiplying the integrand with a $C^\infty$ window centered around $1/2$ and then integrating the nonzero part.

Of course, the remaining integral is still highly oscillatory. One can shrink the support of the window function for increasing $k$, at a rate that depends on the nature of the contributing point. Loosely speaking, one wants to maintain a fixed number of oscillations within the support of the window, such that the integral can be evaluated at a fixed cost for increasing $k$. The size of the support of the window function around a stationary point is typically scaled as $\MO(k^{-1/2})$, whereas for an endpoint contribution it can be smaller: $\MO(k^{-1})$. The size becomes larger for higher-order stationary points, where higher order derivatives of $g$ vanish as well. The more derivatives of $g$ vanish, the less oscillatory the integrand and the larger the contribution. A stationary point $p$ has order $r$ if
\[
 g'(p) = g''(p) = \ldots = g^{(r)}(p) = 0
\]
but $g^{(r+1)}(p) \neq 0$. The associated contribution to the integral has size $\MO(k^{-1/(r+1)})$ and the size of the window scales accordingly. For endpoints and singularities of the integrand, we have $r=0$ and for non-degenerate stationary points, $r=1$.

One typically uses $C^\infty$ window functions in order to avoid introducing discontinuities in the derivatives of the integrand, as these would show up in the asymptotics of the modified integral. An extreme example of a non-smooth window is a block window on $[\lambda, \rho]$, which is discontinuous. This is equivalent to restricting an integral on $[a,b]$ to an integral on $[\lambda,\rho]$. Hence, one expects $\MO(k^{-1})$ spurious contributions from the endpoints $\lambda$ and $\rho$. It is important to note that the windows need not be symmetric.

\section{Asymptotic compression of BEM matrices} \label{Smethod}

We incorporate smooth window functions into the Green's function of an integral equation in order to isolate the contributions of the oscillatory integral operator. The most important consideration is where to place the windows, but first we introduce the specifics of the scheme.

\subsection{Asymptotic compression method}

We consider a domain whose boundary $\Gamma$ is parameterized by $\kappa : \MK \to \Gamma$, such that the single layer potential becomes
\begin{equation*}
  \int_{\Gamma} K(\mathbf{x},\mathbf{y}) v(\mathbf{y}) {\rm d}s(\mathbf{y}) = \int_\MK K(\mathbf{x}, \kappa(\tau)) \Vert \nabla \kappa(\tau) \Vert v(\tau) {\rm d} \tau.
\end{equation*}
In a slight abuse of notation, we identify $v(\mathbf{y})$ with $v(\tau)$. In the following, we shall also use a Green's function in the parameter domain as defined by
\[
 K_\kappa(t,\tau) = K(\kappa(t),\kappa(\tau)) \Vert\nabla\kappa(\tau) \Vert.
\]
For the discretization, we 
use a basis $\Phi_N := \{ \varphi_j \}_{j=1}^N$, defined in the parameter domain $\MK$, of compactly supported functions $\varphi_j$ with support $S_j$. We intend to use piecewise polynomials. Hence,
\begin{equation}
\label{eq:v_N}
 v(\tau) \approx v_N(\tau) = \sum_{j=1}^N c_j \varphi_j(\tau).
\end{equation}
This setting can easily be modified to one where $\Gamma$ itself is approximated by elements, as is common in BEM.

With collocation at the points $x_i = \kappa(t_i)$, $i=1,\ldots,N$, integral equation \eqref{Eie} leads to a discretization matrix $A$ with entries
\begin{equation}
 A_{i,j} = \int_{S_j} K_\kappa(t_i,\tau) \varphi_j(\tau) {\rm d} \tau. \label{EmatEntries}
\end{equation}
We continue with the collocation discretization -- Galerkin is very similar. Recall that a row of $A$ times a column vector $c$ represents an oscillatory integration. Indeed, if $Ac=b$, then
\begin{equation} \label{E:bi}
b_i = \sum_{j=1}^N A_{i,j} c_j = \sum_{j=1}^N c_j \int_{S_j} K_\kappa(t_i,\tau) \varphi_j(\tau) {\rm d} \tau =  \int_{\Gamma} K_\kappa(t_i,\tau) \left(\sum_{j=1}^N c_j \varphi_j(\tau) \right) {\rm d} \tau.
\end{equation}
These integrals only have local contributions near certain critical points, which still depend on the collocation point $t_i$. Assuming we know these, we make this explicit using a smooth bivariate windowing function $w$,
\[
 \tilde{K}_\kappa(t,\tau) = K_\kappa(t,\tau) w(t,\tau).
\]
We introduce some more notation to specify $w$ further. We divide the domain of $w$ in three disjoint sets, $\MK \times \MK = \MK_1 \cup \MK_2 \cup \MK_3$, such that
\[
 w(t,\tau) = \left\{ \begin{array}{ll}
                     1, & \mbox{if~}(t,\tau) \in \MK_1,\\
                     0, & \mbox{if~}(t,\tau) \in \MK_3,
                    \end{array}
\right.
\]
while $w(t,\tau)$ is $C^\infty$ both in $t$ and in $\tau$. 

Denote the compressed matrix corresponding to $\tilde{K}_\kappa$ by $\tilde{A}$, with entries
\begin{equation}
 \tilde{A}_{i,j} = \int_{S_j} \tilde{K}_\kappa(t_i,\tau) \varphi_j(\tau) {\rm d} \tau = \int_{S_j} K_\kappa(t_i,\tau) w(t_i,\tau) \varphi_j(\tau) {\rm d} \tau. \label{EA2entries}
\end{equation}
It is clear that $\tilde{A}_{i,j} = A_{i,j}$ whenever $t_i \times S_j \in \MK_1$, while $\tilde{A}_{i,j} = 0$ if $t_i \times S_j \in \MK_3$. Recall that oscillatory integrals have only local contributions, so the set $\MK_3$ on which the window vanishes will be large. This results in many zeros in $\tilde{A}$, hence \emph{asymptotic compression}.

The matrix vector product $\tilde{A}c$ corresponds to the following oscillatory integral:
\begin{equation}
\label{eq:discrete_oscillatory_integral}
	\sum_{j} \tilde{A}_{i,j} c_j = \int_\MK K_\kappa(t_i,\tau) w(t_i,\tau) \sum_{j} c_j \varphi_j(\tau) {\rm d} \tau = \int_\MK K_\kappa(t_i,\tau) w(t_i,\tau) v_N(\tau) {\rm d} \tau.
\end{equation}
The two integrals, $\tilde{A}c$ above and $b_i$ in \cref{E:bi}, do not agree in general. However, the assumption is that when $v_N$ is an oscillatory function, and since the Green's function is oscillatory as well, for well chosen support sets $\MK_1$ and $\MK_3$, the two integrals do agree asymptotically. Their difference is an asymptotic error which is spectrally small in the wavenumber $k$ even though we consider large perturbations of the matrix $A$, which may have a large condition number.

We proceed by replacing $A$ by $\tilde{A}$ in \eqref{Esystem}. This leads to a new linear system of equations,
\[
 \tilde{A} \tilde{c} = b,
\]
with a sparse matrix $\tilde{A}$. For collocation, we simply have $b_i = -u^{\text{inc}}(t_i)$. In this case, the weight function $w(t,\tau)$ need only be defined for the discrete points $t_i$.

\subsection{Explicit expression of a smooth windowing function} \label{SchoiceWind} 

Determining suitable locations for the window functions is a crucial problem, due to the fact that the phase of the solution is not known. We discuss this separately in the next section. Here, we describe the mathematical expression of the function.

The window function is a function of two variables, $t$ and $\tau$, but the integration is performed with respect to $\tau$. Therefore, we consider the location of the window to be a function of $t$. In general, $w(t,\tau)$ is a sum of elementary window functions with locations depending on $t$:
\begin{equation*}
 w(t,\tau) = \sum_l \chi(\tau,\lambda_l(t), l_l(t), r_l(t), \rho_l(t)).
\end{equation*}
They are supported on $[\lambda,\rho]$ and equal to $1$ on the smaller interval $[l,r]$. When two subsequent windows $l$ and $l+1$ in $w$ overlap or come too close, we join them into $\chi(\tau,\lambda_l, l_l, r_{l+1}, \rho_{l+1})$. The specific choice of $\chi$ can be important from a numerical point of view, but we expect very similar results as long as all derivatives have a reasonable bound. The extreme case would be a block window ($\lambda=l$ and $r=\rho$), introducing the aforementioned spurious endpoint contributions. We choose $\chi$ as the elementary window also used in \cite{bruno}, and shown in \cref{FoscIntIll}:
\begin{equation}
 \chi(\tau, \lambda, l, r, \rho) = \begin{cases}  0, & \tau \leq \lambda \\ \exp\left(\frac{2e^{\frac{l-\lambda}{\tau-l}}}{\frac{\tau-l}{\lambda-l}-1}\right), & \tau \in (\lambda,l) \\ 1, & \tau \in [l,r]  \\ \exp\left(\frac{2e^{\frac{r-\rho}{\tau-r}}}{\frac{\tau-r}{\rho-r}-1}\right), & \tau \in (r,\rho) \\ 0, & \tau \geq \rho. \end{cases} \label{Ewin2D} 
\end{equation}

The width of the window will be determined automatically further on, even taking into account possible asymmetry. One expects the window to optimally scale as mentioned in the previous section on the asymptotics of oscillatory integrals, i.e. as ${\mathcal O}(k^{-1/(r+1)})$ around a stationary point of order $r$. Single reflections of rays bouncing off an obstacle are associated with stationary points of order $1$ \cite{groth2017highorder}. Degenerate stationary points do also arise frequently in scattering problems, for example at shadow boundaries where an incoming ray is tangential to $\Gamma$. Grazing rays, which give rise to creeping rays travelling along the boundary of the obstacle, are associated with stationary points of order $2$ in asymptotic analysis \cite{unifPhaseExtr,MelroseTaylor}. Note that these points depend on the angle of the incoming wave. Hence, so do the contributing points and so does our computational scheme.

\subsection{Robustness for low-order discretizations and computation of matrix entries}

The asymptotic analysis of integrals of the form \eqref{eq:oscillatoryintegral} usually assumes a smooth integrand. Typically, later terms in the asymptotic expansion depend on higher order derivatives of the integrand. Yet, in our setting we are faced with the discretized integral \eqref{eq:discrete_oscillatory_integral} that is not particularly smooth: the smoothness is limited at least by the smoothness of the basis functions $\varphi_j$. We can only expect low-order asymptotics because of the piecewise polynomial basis functions used in this article. In this sense, it is somewhat unexpected that the results further on show that asymptotic compression is also effective for piecewise linear basis functions, as smoothness of the integrand is important for asymptotic behaviour as $k\rightarrow \infty$.

This apparent robustness also allows a significant simplification, that makes the implementation much cheaper. Rather than integrating the window function exactly, we will explore approximating the compressed matrix element of $\tilde{A}$ simply by rescaling the uncompressed element of $A$:
\begin{align}
 \tilde{A}_{i,j} &= \int_{S_j} K_\kappa(t_i,\tau) w(t_i,\tau) \varphi_j(\tau) {\rm d} \tau \nonumber  \\
 &\approx w(t_i, \tau_j) \int_{S_j} K_\kappa(t_i,\tau) \varphi_j(\tau) {\rm d}\tau \nonumber \\ 
 &= w(t_i, \tau_j) A_{i,j}. \label{eq:element_approximation} 
\end{align}
Here, $\tau_j$ is a point within the support of $\varphi_j$. This is a rather crude approximation of the exact weighted integral. However, it has the important advantage that the element of $\tilde{A}$ is a simple rescaling of the corresponding element of $A$. Thus, existing code for the computation of entries of $A$ can maximally be reused.

Note that this approach is still decidedly different from simply compressing $A$ by putting entries to zero. Even with our simplified approach, the elements of $A$ are weighted, and thus decay to zero with some regularity. Simply discarding some matrix entries would correspond to using block-window functions, which \cref{SaddRes} indicates still works but with much larger errors.

\section{Adaptive location of the window functions}
\label{s:adaptive}

We consider three ways of finding a suitable support of the window functions. The first is based on the simple geometric notion of visibility. This approach is inferior to the following ones, but it does illustrate that determining supports is much easier than performing ray tracing. Second, we numerically estimate suitable windows by solving a complete scattering problem at a moderate value of the wavenumber. These windows are kept constant and reused for larger frequencies. Finally, we compute windows adaptively in a frequency sweep, refining the windows for increasing frequency. The latter scheme is most efficient, and we can start at a smaller value of the wavenumber for the initial windows.

\subsection{Visibility criterion} \label{SvisExplain}
Asymptotically, a wave field can be seen as a collection of rays along straight paths (in a homogeneous medium)~\cite{babich1991shortwave}. Points of reflections of these rays physically correspond to stationary points in the integral equation formulation. For scatterers with a moderately complicated geometry, it becomes unwieldy to trace all reflecting rays through the scene. However, in the integral equation formulation, the oscillatory integral in equation \eqref{Eie} corresponding to a point $\mathbf{x} \in \Gamma$ on the boundary always asymptotically corresponds to a sum of contributions from a limited number of points.

If a point $B$ on the boundary is not directly visible from a point $A$ on the illuminated part of the boundary, then $B$ can not be a stationary point of the oscillatory integral at $A$. Indeed, there can be no \emph{direct} ray from $B$ to $A$ (or vice-versa). There can only be indirect connections via intermediate reflection points $C$ elsewhere on the boundary, or via creeping rays along the surface of the obstacle. In the former case, $B$ is a stationary point for the integral at $C$ which is subsequently linked with $A$, but there is no direct connection from $B$ to $A$ via a stationary point which has to be taken into account. This means that for any given point $A$ on the illuminated part of the boundary, we can asymptotically discard the density in all non-visible parts of the boundary from that point $A$, unless that part is close to $A$ due to the contribution from the singularity.

Note that this simple criterion does not require explicit information on the phase of the density $v(\mathbf{y})$. However, this reasoning does not immediately hold for points on the shadow side of the boundary, for which the situation is more complicated (see \S\ref{Sshadow} further on).

In a collocation solver, we have to decide for each collocation point $t_i$ where to place the window function $w(t_i,\tau)$. For objects with self-reflections, one could trace a ray from the source through the scattering scene to the collocation point, much like ray tracing in computer graphics. Each point of reflection is expected to be a stationary point for the next point of reflection, or ultimately for the collocation point. Hence, each of these points needs a window around it, and any two subsequent points have to be able to `see' each other for a reflection to be possible. We do not consider any ray tracing in our visibility criterion, but simply exploit the latter visibility condition. Computing which points can see one another gives the regions where critical points can lie. The visibility problem is also common in computer graphics on the GPU, and the computation can be sped up using Kd-trees or Bounding Volume Hierarchies \cite{starRaytr}.

\subsection{Adaptive asymptotic compression} \label{Scorr}
Points that contribute to a general oscillatory integral can be detected automatically using a sliding window $\zeta(x-z)$ centered around a point $z$. If the window contains no contributing points, asymptotic theory predicts the windowed integral
\begin{equation}
\label{Ecorrelation}
 F(z) = \int_a^b \zeta(x-z) f(x) e^{i k g(x)} {\rm d}x 
\end{equation}
to be superalgebraically small in the frequency parameter $k$ (for $C^\infty$ integrands).
As explained in \cref{s:integrals}, the integrand would look like the red dashed line at the right of \cref{FoscIntIll} where the positive and negative parts cancel. Conversely, the function $F(z)$ can be expected to be large if the window around $z$ does contain a contributing point. This corresponds to shifting the solid green window function to $z$ near $1/2$ in \cref{FoscIntIll}, where the integrand is locally less oscillatory.

In the context of integral equation methods, one can compute correlations between any two points on the boundary as integrals of the unmodified integrand times a moving window function. This can only be done a posteriori, since the integrand in~\eqref{Eie} involves the density function which is unknown a priori. 

\begin{figure}
\centering
\begin{tikzpicture}[scale=1.6]
	\draw (0,0) node {$\begin{pmatrix} A_{1,1} & A_{1,2} & \cdots & \cdots & \cdots & A_{1,N} \\ A_{2,1} & A_{2,2} &  & \cdots & & A_{2,N} \\ \vdots & \vdots &  & \ddots & & \vdots \\ A_{N,1} & A_{N,2} & & \cdots & & A_{N,N} \end{pmatrix}$};
	\draw (1.6,0) node {$\begin{pmatrix} c_1 \\ c_2 \\ \vdots \\ c_N \end{pmatrix}$};
	\draw (3.8,0) node (R) {$= \begin{pmatrix} R_{1,1} & \cdots & {\color{purple} R_{1,q}} & {\color{orange} R_{1,q+1}} & \cdots &  R_{1,Q} \\ R_{2,1} & \ddots & &  & & R_{2,Q} \\ \vdots & & & & \ddots & \vdots \\ R_{N,1} & &  & & \cdots & R_{N,Q} \end{pmatrix}$}; 
	\coordinate (r11) at (3.35, 0.6); 
	\draw (3.5, 0.8) node {$\leftarrow \sigma_q \rightarrow$};
	\draw (5.6, -0.1) node {\begin{minipage}{0.3cm} \begin{align*} & \uparrow \\ & t_n \\ & \downarrow \end{align*} \end{minipage} };
	\draw (-1.5, -0.1) node {\begin{minipage}{0.3cm} \begin{align*} & \uparrow \\ & t_n \\ & \downarrow \end{align*} \end{minipage}};
	\coordinate (ro) at (3.6, 0.6);
	\draw[purple,rounded corners] (-1.2, 0.3) rectangle (0.3, 0.5) node[yshift=-3] (ap) {};
	\draw[orange,rounded corners] (-0.7, 0.28) rectangle (0.5, 0.52) node (ao) {};
	\draw[purple,rounded corners] (1.45, -0.11) rectangle (1.75, 0.5) node (cp) {};
	\draw[orange,rounded corners] (1.43, -0.2) rectangle (1.77, 0.25) node (co) {};

	\path[purple,->] (ap) edge[out=5, in=135,middlearrow={0.3}{latex}] (r11);
	\path[purple,->] (cp) edge [out=15, in=135,middlearrow={0.5}{latex}] (r11);

	\draw[purple, very thick, domain=0.001:0.2499, samples=100]  plot ({\x*5-1.2}, {0.6+0.5*exp(2*exp(0.25/(\x-0.25))/((\x-0.25)/(-0.25)-1) )}); 	
	\draw (0.05, 1) -- (0.05, 1.2) node[above,purple] {$\sigma_q$};
	\draw[purple, very thick, dashed, domain=0.7501:0.999, samples=100]  plot ({\x-0.7}, {0.6+0.5*exp(2*exp(-0.25/(\x-0.75))/((\x-0.75)/(0.25)-1) )});  
	\path[purple,->] (0.2,0.85) edge [out=0, in=135,middlearrow={0.35}{latex}] (r11);

	\draw[orange, very thick, domain=0.001:0.2499, samples=100]  plot ({\x*4.4-0.7}, {0.6+0.5*exp(2*exp(0.25/(\x-0.25))/((\x-0.25)/(-0.25)-1) )});
	\draw (0.4, 1) -- (0.4, 1.2) node[above right,orange] {$\sigma_{q+1}$};
	\draw[orange, very thick, dashed, domain=0.7501:0.999, samples=100]  plot ({\x*0.4+0.1}, {0.6+0.5*exp(2*exp(-0.25/(\x-0.75))/((\x-0.75)/(0.25)-1) )}); 
	
	\draw[->] (-1.2, 0.5) -- (-1.2, 1.3) node[left] {$\zeta$}; 
	\draw[->] (-1.3, 0.6) -- (1.2, 0.6) node[above left] {$t_l$}; 
\end{tikzpicture}
\caption{Schematic depiction of the computation of the correlations $R_{n,q}$ given by \cref{EroisSum}. The matrix and the solution vector $c$ %PO Matrix $A$ and solution vector $c$ 
are computed for a wavenumber $k_1$. Each row $n$ of $A$ is multiplied by a sliding window centered around $\sigma_q$, and the inner product with the solution vector $c$ results in a localised correlation value $R_{n,q}$. Recall that row $n$ correspond to collocation point $t_n$. For more implementation details, we refer the reader to the code available online \cite{github}.}
\label{FcorrProc}
\end{figure}
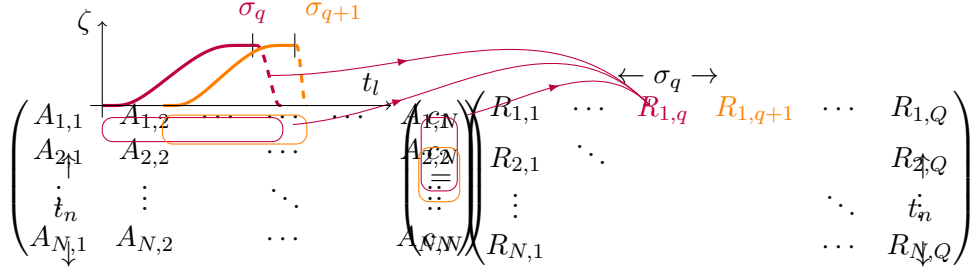

Having computed $v_N$ for a (moderate) value of the wavenumber $k_1$, we can compute correlations comparable to \eqref{Ecorrelation} above. We do so for each collocation point $t_n^{(k_1)}$, $m=1,\ldots,N$, and use windows around center points $z=\sigma_q$, $q=1, \ldots, Q$ where we choose $Q=1.5N$. We use more centers than collocation points by a factor of $1.5$ in order to determine more accurately where the important contributions are located. The window function we use around the center $\sigma_q$ is 
\[
\zeta(\tau-\sigma_q) = \chi(\tau-\sigma_q, -T,0, 0, T).
\]
We compute a correlation matrix $R$, with elements given by:
\begin{align}
	R_{n,q} & = \int_\MK \zeta(\tau-\sigma_q)  K_\MK(t_n, \tau) v_N(\tau) {\rm d}\tau \\ 
	& = \sum_{j=1}^N c_j \int_\MK \zeta(\tau-\sigma_q)  K_\MK(t_n, \tau) \varphi_j(\tau) {\rm d}\tau \nonumber \\
	& \approx \sum_{l=1}^N \zeta(t_l - \sigma_q) A_{m,l} c_l. \label{EroisSum}
\end{align}
Recall from \eqref{eq:v_N} that $v_N(\tau)$ is our approximation to the density function. In the latter step, we have again approximated a windowed integral by a weighted sum of matrix elements of the dense discretization matrix $A$, which makes the computation of all correlations a cheap operation. No additional numerical integration is required beyond what was already computed for the regular BEM matrix $A$. This computation is visualised at the top of \cref{FcorrProc}: as $\sigma_q$ shifts, other parts of the matrix-vector product $Ac$ will be windowed.

For the compression phase at a higher wavenumber $k_2 > k_1$, we identify the collocation point $t_i^{(k_2)}$ with the closest collocation point $t_n^{(k_1)}$ at the lower wavenumber $k_1$ at which the correlations were computed and we identify the interaction point $t_j^{(k_2)}$ with the closest window center $\sigma_q^{(k_1)}$. We retain all unmodified matrix elements corresponding to correlations higher than a fixed percentage $\xi$ times the maximal correlation on each row. It is possible to determine suitable threshold percentages $\xi$ adaptive to the geometry by monitoring the relative compression error $\Vert \tilde{A}c-b\Vert / \Vert b\Vert$, or even only for a few representative elements of $b$. One can also compare correlations to $\xi$ times the maximal $|R_{n,q}^{(k_1)}|$ for all $m$ and $n$, but we found a row-wise threshold to be more robust with respect to changes in the maximum between rows. 

We devise a window that covers the retained centers, and smoothly decays to zero away from them. The decay given by the weight function $w(t_i,\tau)$ in \cref{EA2entries} is defined by finding a $t$ such that $|\sigma_q - \sigma_t| < T$ and $|R_{n,t}| \geq \xi \max_l |R_{n,l}|$ %PO $|R_{1,t}| \geq \xi \max_q |R_{1,q}|$
as depicted in \cref{FA2Proc}. If no such $t$ exists, the compressed matrix element $\tilde{A}_{i,j}$ is zero, thus introducing sparsity. We also always include the Green's function singularity $t_i$, even when the computed correlations are small there, since otherwise $\tilde{A}$ may not always have full rank. Care is taken in the implementation to merge windows that would otherwise overlap.

\begin{figure}
\centering
\begin{tikzpicture}[scale=1.6]
	\draw (2,-1.3) node {$\tilde{A}_{{\color{brown} i}, {\color{green} j} }^{(k_2)} =  \begin{cases} A_{{\color{brown} i}, {\color{green} j} }^{(k_2)} & \text{if } |{\color{purple} R}_{{\color{brown} n}, {\color{green} q}}| \geq \xi \displaystyle \max_l |{\color{purple} R}_{{\color{brown} n},l}|  \\ \zeta({\color{green} t_j} - {\color{brown} t_i}) A_{{\color{brown} i}, {\color{green} j} }^{(k_2)} & \text{if } |{\color{green} t_j} -{\color{brown} t_i}| < T \\ \zeta({\color{green} t_j} - {\color{blue} \sigma_t}) A_{{\color{brown} i}, {\color{green} j} }^{(k_2)} & \text{if } 0 < |{\color{green} \sigma_q} - {\color{blue} \sigma_t}| < T \\ 0 & \text{otherwise.}  \end{cases}$};  %PO \draw (2,-1.3) node {$\tilde{A}_{{\color{brown} i}, {\color{green} j} }^{(k_2)} =  \begin{cases} A_{{\color{brown} i}, {\color{green} j} }^{(k_2)} & \text{if } |{\color{purple} R}_{{\color{brown} n}, {\color{green} q}}| \geq \xi \displaystyle \max_l |{\color{purple} R}_{{\color{brown} n},l}|  \\ \zeta({\color{green} t_j} - {\color{blue} \sigma_t}) A_{{\color{brown} i}, {\color{green} j} }^{(k_2)} & \text{if } \displaystyle \min_{\displaystyle {\color{blue} t} \text{   s.t. } \displaystyle |{\color{purple} R}_{{\color{brown} n},{\color{blue} t} }| \geq \xi \displaystyle \max_l |{\color{purple} R}_{{\color{brown} n},l}| \text{ or } {\color{blue} \sigma_t} \approx {\color{brown} t_i}} |{\color{green} \sigma_q} - {\color{blue} \sigma_t}| < T   \\ 0 & \text{otherwise.}  \end{cases}$}; 
	%PO \path[blue,->] (0.2,-1.15) edge [out=90, in=90,middlearrow={0.5}{latex},looseness = 0.4] (4.7,-1.15);

	\draw[brown] (-0.4,-2.3) node (tmti) {$\displaystyle \min_n |t_n^{(k_1)} - t_i^{(k_2)}|$};
%PO changed coordinates and added paths from here
	\path[brown,->] (0.05,-1.45) edge [out=270, in=90,middlearrow={0.5}{latex},middlearrow={1.0}{latex}] (tmti);
	\path[brown] (2.65,-0.87) edge [out=90, in=45,middlearrow={0.5}{latex},looseness = 0.6] (-0.2,-1);
	\path[brown,->] (-0.2,-1) edge [out=225, in=90,looseness = 1] (tmti);%\path[brown,->] (0.05,-1.45) edge [out=270, in=130,middlearrow={1.0}{angle 90}] (tmti);%	\path[brown,->] (0.05,-1.45) edge [out=270, in=130,middlearrow={0.5}{latex},middlearrow={1.0}{latex}] (tmti);%	\path[brown,->] (2.6,-0.87) edge [out=90, in=130,middlearrow={0.35}{latex},looseness = 0.8] (tmti);
	%\path[brown,->] (3,-1.2) edge [out=160, in=90,middlearrow={0.5}{latex}] (tmti);

	\draw[green] (1.8,-2.3) node (signtj) {$\displaystyle \min_q |\sigma_q^{(k_1)} - t_j^{(k_2)}|$};
	\path[green] (0.2,-1.45) edge (0.2,-1.65);
	\path[green,->] (0.2,-1.65) edge [out=270, in=90,middlearrow={0.25}{latex},looseness = 0.75] (signtj);
	\path[green,->] (1.05,-1.6) edge [out=270, in=90,middlearrow={0.35}{latex},looseness = 1] (signtj); %\path[green,->] (1.05,-1.35) edge [out=270, in=90,middlearrow={0.55}{latex},looseness = 0.6] (signtj);
	\path[green,->] (5.3,-2.17) edge [out=90, in=90,middlearrow={0.7}{latex},looseness = 0.18] (signtj);

	\draw[blue] (4.5,-2.3) node (sigtq) {$\displaystyle \min_{\displaystyle {\color{blue} t} {\color{black} \text{   s.t. }} \displaystyle {\color{black} |}{\color{purple} R}_{{\color{brown} n} {\color{black} ,} {\color{blue} t} } {\color{black} | \geq \xi \displaystyle \max_l |} {\color{purple} R}_{{\color{brown} n} {\color{black} ,l} }{\color{black} |} } {\color{black} |} {\color{green} \sigma_q} {\color{black} -} {\color{blue} \sigma_t}{\color{black} |}$}; %\draw (4.5,-2.3) node (sigtq) {$\displaystyle {\color{blue} \min}_{\displaystyle {\color{blue} t} \text{   s.t. } \displaystyle |{\color{purple} R}_{{\color{brown} n},{\color{blue} t} }| \geq \xi \displaystyle \max_l |{\color{purple} R}_{{\color{brown} n},l}| } |{\color{green} \sigma_q} - {\color{blue} \sigma_t}|$}; 
	\path[blue,->] (1.45,-1.6) edge [out=270, in=90,middlearrow={0.25}{latex},looseness = 0.4] ([xshift=-0.4cm,yshift=0.2cm]sigtq); %(sigtq) ++ (-1.0,0); %([xshift=-0.5]sigtq);
	\path[blue,->] (3.33,-1.6) edge [out=270, in=90,middlearrow={0.35}{latex},looseness = 0.6] ([xshift=-0.4cm,yshift=0.2cm]sigtq);
\end{tikzpicture}
\caption{The computation of the compressed matrix $\tilde{A}^{(k_2)}$ at wavenumber $k_2$, with entries given by \cref{EA2entries} and approximated by \cref{eq:element_approximation}. The correlations $R_{n,q}$ were computed at a smaller wavenumber $k_1$, see \cref{FcorrProc}. Suitable values of $n$ and $q$ to use for $\tilde{A}^{(k_2)}_{i,j}$ are found by matching the collocation point $t_i^{(k_2)}$ at wavenumber $k=k_2$ to collocation points and window centers at $k=k_1$. No compression is applied if the correlation is large. The matrix element is put to zero where the correlation is low, and the element is weighted in between.}
\label{FA2Proc}
\end{figure}

The goal is to maintain a compression error of the same order as the discretization error. Without any assumption on the density $v_N(\tau)$, one might expect large errors from throwing away a part of the matrix. However, if the incident wave changes to give significant contributions elsewhere, our correlations $R_{m,n}$ will increase there as we compute the density at a moderate frequency. This reasoning may lead to a rigorous proof that our method only removes negligible quantities, and we provide numerical results further on for various types of incident waves.

Note that this procedure leads to a window for each collocation point at the wavenumber $k_1$ where the correlations are computed. At larger frequencies, the discretization may be refined and there may be more collocation points. For each of those points, we simply use the window of the closest collocation point at $k_1$. We can reuse the correlation matrix $R_{m,n}$ for larger $k$ because the list of stationary points or contributing regions that $R_{m,n}$ provides has a structure which is essentially independent of $k$. That is because stationary points are essentially a geometrical feature, and hence their location is independent of the frequency.

Note that the approach using correlations has an additional advantage. When determining the asymptotic expansion of the oscillatory integral as in \cite{MelroseTaylor}, one computes the points where the gradient of the phase of the integrand is exactly zero and determines the contributions coming from them. This is implicitly done by the reflecting points in ray tracing. The correlations also find regions where the gradient of the phase is small (in absolute value) but does not cross zero. These regions have no contribution as $k\rightarrow\infty$, but a finite $k$ might have to be very large for this contribution to become negligible. Thus, including these regions is also advantageous, which the correlation technique does. This has also been remarked in \cite[\S 5.1.1]{ganesh}. We also retain contributions from creeping rays if they are still significant at the frequency at which the correlations are computed.

We illustrate the method in \cref{Snum} with numerical results, but we already point out two issues. The first issue relates to the value of the `moderate' wavenumber $k_1$. The adaptive asymptotic compression method is adaptive with respect to the geometry and computes the correlation matrix $R$ only at this lowest wavenumber. These correlations provide valuable high-frequency information only as long as there is `enough' high-frequency behaviour present at that value of $k$. Solving the full boundary element problem at this value of $k$ might thus be expensive as well. Indeed, this is the problem we are trying to speed up in the first place. The second issue is that thus far we are using fixed windows for each value of $k$, computed at the small value of $k_1$. In doing so, we are not exploiting the fact that windows could be smaller at higher frequencies as described in \cref{SchoiceWind}.

Both of these issues are adressed by adaptively recompressing the matrix as the frequency increases. This is especially appealing when one has to solve the same problem for a range of frequencies. As will be explained in \cref{SadapRecomprExplain}, we can take the smallest frequency to be smaller, and we can efficiently obtain smaller windows for larger values of $k$.

\subsection{Adaptive recompression at higher frequencies} \label{SadapRecomprExplain}

In the adaptive asymptotic compression technique thus far, we have been using fixed windows computed at a moderate wavenumber. This leads to a fixed compression ratio, independent of the frequency and a decreasing asymptotic error. However, it should be possible to use smaller windows at larger frequencies for a constant error, hence improving compression with increasing frequency. In the absence of prior asymptotic information, one way to achieve higher compression at higher frequencies in general is to recompute the correlations $R_{m,n}$ at each wavenumber.

Of course, computing the full correlation matrix is prohibitively expensive, since it requires the dense matrix $A$. However, if we have already computed a compressed matrix $\tilde{A}^{(k_j)}$ %PO $\tilde{A}^{(j)}$
at a smaller wavenumber $k_j$, then we can omit a large part of the computation. Indeed, if the correlation between two points on the boundary is small at a certain frequency, indicating that they are not asymptotically linked, it is safe to assume that their correlation will only be smaller at higher frequencies. Thus, we have to compute correlations only for the non-zero entries of $\tilde{A}^{(k_j)}$. Assuming some of these will be small, because the integrals are more oscillatory, this leads to a sparser matrix $\tilde{A}^{(k_{j+1})}$ at a higher wavenumber $k_{j+1}$. We can repeat this reasoning, where at each higher frequency we benefit from the improved sparsity of the previous stage.

To be precise, we introduce some notation. Say we have a range of strictly increasing wavenumbers $k_j$, $j=1,\ldots,J$, and we use $N_j$ degrees of freedom at each wavenumber $k_j$. In a full correlation technique, we would use $1.5 N_j$ centers $\sigma_{q}^{(k_j)}$, $q=1,\ldots,Q_j$, for each collocation point $t_{n}^{k_j}$, $n=1,\ldots,N_j$. The corresponding correlation matrix $R_{n,q}^{(k_j)}$ determines window functions $w_{n}^{(k_j)}$ to compute the compressed matrix $\tilde{A}^{(k_j)}$, and the solution of the discretized integral equation is $\tilde{c}^{(k_j)} = \tilde{A}^{(k_j)} \setminus b^{(k_j)}$. %$\tilde{c}^{(j)} = \tilde{A}^{(j)} \setminus b^{(j)}$.

In the adaptive asymptotic recompression scheme, we only compute the full correlation matrix $R^{(k_1)}$ at the smallest wavenumber $k_1$. We reuse $R^{(k_1)}$ for the next few wavenumbers $k_l$ until and including $k_j=2k_1$. For this larger wavenumber $k_j$, after having computed the solution vector $\tilde{c}^{(k_j)}$ and the associated density function $v^{({N_j})}$, we only compute some of the elements of $R^{(k_j)}$: %PO$R^{(j)}$: 
\[
 R_{n,q}^{(k_j)} = \left\{ \begin{array}{cl}
                \int_\MK \zeta\left(\tau-\sigma_{q}^{(k_j)}\right)  K_\MK(t_n^{(k_j)}, \tau) v^{(N_j)}(\tau) {\rm d}\tau \nonumber,  \quad &\textrm{if~} \sigma_q^{(k_j)} \in \SUPP w_{m'}^{(j-1)}, \\
                0, & \textrm{otherwise}.
               \end{array}\right.
\]%PO  R_{j,m,n} = \left\{ \begin{array}{cl}
where $m'$ is such that $t_{m'}^{(k_{j-1})}$ %$t_{m'}^{(j-1)}$ 
is the closest collocation point to $t_{m}^{(k_j)}$. Subsequently, the element $R_{m,n}^{(k_j)}$ can be approximated by a weighted sum of elements of $\tilde{A}^{(k_j)}$ instead of $A^{(k_j)}$ as before in \cref{EroisSum}. Finally, new windows are computed based on thresholding the correlations using a percentage $\xi$, and these are used for the compressed matrices $\tilde{A}^{(k_l)}$ at the next wavenumbers $k_l > k_j$ until the wavenumber has doubled again, so $k_l=2k_j$ and so on. This loop can be visualised by inserting $R_{n,q}$ of \cref{FcorrProc} into \cref{FA2Proc} and by inserting $\tilde{A}^{(k_l)}$ of the latter as $A^{(k_j)}$ into the former for $k_l=2k_j$. 

One may believe that the cost of repeatedly computing correlations remains substantial. However, we limit ourselves to recomputing correlations at wavenumbers that are logarithmically spaced. As a result, the cost of computing correlations as described above is swamped by the cost of solving multiple problems when traversing a frequency range as numerical results show later on. Moreover, we found it to be very robust in the way it is described above. Simplifications in which we count the number of nearby zero crossings or in which we refine our windows near the maximum of correlations, assuming this maximum corresponds closely to a stationary point, are cheaper but less accurate.

\subsection{Collocation points in the shadow region(s) and Fermat's principle} \label{Sshadow}

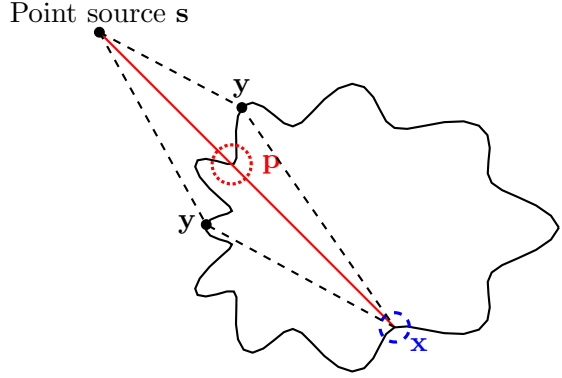
\begin{figure}[t]
\centering
\begin{tikzpicture}[scale=1.3]
 \draw[fill=black] (-2,2) circle[radius=0.05] node[right,above] {Point source $\mathbf{s}$};
 \draw [thick,  domain=0:360, samples=100]  plot ({(0.7+0.01*abs(\x-180)+0.2*cos(10*\x))*cos(\x)}, {(0.7+0.006*abs(\x-180)+0.2*cos(10*\x))*sin(\x)} );
 \coordinate (x) at ({(0.7+0.01*abs(305.528-180)+0.2*cos(10*305.528) )*cos(305.528)}, {(0.7+0.006*abs(305.528-180)+0.2*cos(10*305.528))*sin(305.528)} );
 \coordinate (y) at ({(0.7+0.01*abs(128.723-180)+0.2*cos(10*128.723))*cos(128.723)}, {(0.7+0.006*abs(128.723-180)+0.2*cos(10*128.723))*sin(128.723)} );
 \draw[thick,red] (-2,2) -- (y) -- (x);
 \draw[red,very thick, densely dotted] (y) circle[radius=0.2] node[red, right] {~~$\mathbf{p}$};
 \draw[blue,very thick, dashed] (x) circle[radius=0.15] node[blue,below] {~~~~~$\mathbf{x}$};
 \def\phit{178}
 \coordinate (z) at ({(0.7+0.01*abs(\phit-180)+0.2*cos(10*\phit))*cos(\phit)}, {(0.7+0.006*abs(\phit-180)+0.2*cos(10*\phit))*sin(\phit)} );
 \draw[fill=black] (z) circle[radius=0.05] node[left] {$\mathbf{y}$};
 \draw[thick,dashed] (-2,2) -- (z) -- (x);
 \def\th{110};
 \coordinate (q) at ({(0.7+0.01*abs(\th-180)+0.2*cos(10*\th))*cos(\th)}, {(0.7+0.006*abs(\th-180)+0.2*cos(10*\th))*sin(\th)} );
 \draw[fill=black] (q) circle[radius=0.05] node[above] {$\mathbf{y}$};
 \draw[thick,dashed] (-2,2) -- (q) -- (x);
\end{tikzpicture}
\caption{A collocation point in the shadow region ($\mathbf{x}$ in the figure) always has a stationary point $\mathbf{p}$ on the boundary in the illuminated region, that lies on the line connecting $\mathbf{x}$ to the source $\mathbf{s}$. By Fermat's principle this point minimizes the sum $\Vert \mathbf{s}-\mathbf{y} \Vert + \Vert \mathbf{y} - \mathbf{x} \Vert$ for $\mathbf{y}$ on the boundary.}
\label{Fpointsource}
\end{figure}

In some cases, the location of stationary points has been extensively studied, most notably for single convex smooth scatterers in 2D and 3D~\cite{bruno,Dominguez,ganesh,Ganesh2}. In this setting, one knows the phase of the solution asymptotically, and hence one can study the exact location of all stationary points. The number of stationary points depends on the region in which the collocation point lies: the part of the obstacle that is visible by the incoming wave, the part that is in the shadow, or the transitional part in between. It appears that collocation points in a shadow region are always associated with one (or more) stationary points in the illuminated region. These points will be found automatically using the adaptive scheme above, but should be taken into account explicitly when using the visibility criterion.

As an example of the references given above, the computation can be done exactly for an incident plane wave in the $x$ direction on a circle of radius $1$, parametrised as $\kappa(\tau) = [\cos(2\pi\tau), \sin(2\pi\tau)]$ for $\tau \in [0,1]$. Each collocation point has at least one critical point given by $\tau = t$, which corresponds to the singularity of the Green's function. For points in the illuminated region and in the shadow boundary, this is the only one. There is exactly one other critical point, a stationary point, when the collocation point lies in the shadow region, and it arises as follows. One can assume that the phase of $v(\tau)$ on the illuminated side is equal to the phase of the incident wave \cite[Rem. 4.1]{sparseDiscr}. The phase of $v(\tau)$ is thus equal to $\cos(2\pi \tau)$ for $\mathbf{y}=\kappa(\tau)$ in the illuminated region, and the phase of $K_\kappa(t,\tau)$ is the Euclidean distance to the collocation point. Hence, there exists a slowly varying $f$ such that
\[
 K_\kappa(t,\tau)v(\tau) = f(t,\tau,k) \exp[ik g(t,\tau)]
\]
with
\[
 g = \cos2\pi \tau +\sqrt{(\cos 2\pi t-\cos2\pi \tau)^2 + (\sin 2\pi t-\sin2\pi \tau)^2}
\]
and therefore
\[
 \frac{\partial g}{\partial \tau} = \frac{-\sin2\pi \tau}{2\pi} +\frac{\sin2\pi \tau\cos 2\pi t -\cos2\pi \tau\sin 2\pi t}{2\pi \sqrt{(\cos 2\pi t-\cos2\pi \tau)^2 + (\sin 2\pi t-\sin2\pi \tau)^2} }.
\]
The derivative of the phase is zero at $\tau = 1/2-t$, hence $1/2-t$ is a stationary point $\mathbf{p}$ of the integral \eqref{Eslpot}. It lies in the illuminated region on the horizontal line through the collocation point.

While the computation of the circle is specific, the conclusion appears to be generic: it holds for more general objects that for any collocation point in the shadow region, there is a stationary point in the illuminated region on the line connecting the collocation point with the source. In \cref{Fpointsource}, we illustrate the existence of a stationary point on this line according to Fermat's principle, which states that rays travel along paths that minimize the total travel time. Consider a point $\mathbf{x}$ in the shadow region, and a line from that point to the point source (or to infinity in case of a plane wave). Among all possible points on the boundary $\mathbf{y}$, the distance from the point source to $\mathbf{y}$ plus the distance from $\mathbf{y}$ to $\mathbf{x}$ is trivially minimized precisely by the point $\mathbf{p}$ shown, where the line crosses the illuminated region.

Since we have assumed a non-penetrable scatterer in this paper, no ray physically passes through the interior of the obstacle. The stationary point arises because the Green's function only takes into account the distance between two points, and not the relative location of those points to the obstacle. Intuitively, one can say that the density $v(\mathbf{y})$ is asymptotically very small in the shadow region, hence the integral in \cref{Eie} needs a contribution from the illuminated region to satisfy the boundary condition. Our adaptive strategy for general geometries confirms this and automatically recovers the stationary point.

\section{Numerical results} \label{Snum}

\begin{figure}[t]
\centering
\hspace*{-2cm}\includegraphics[width=1.3\linewidth]{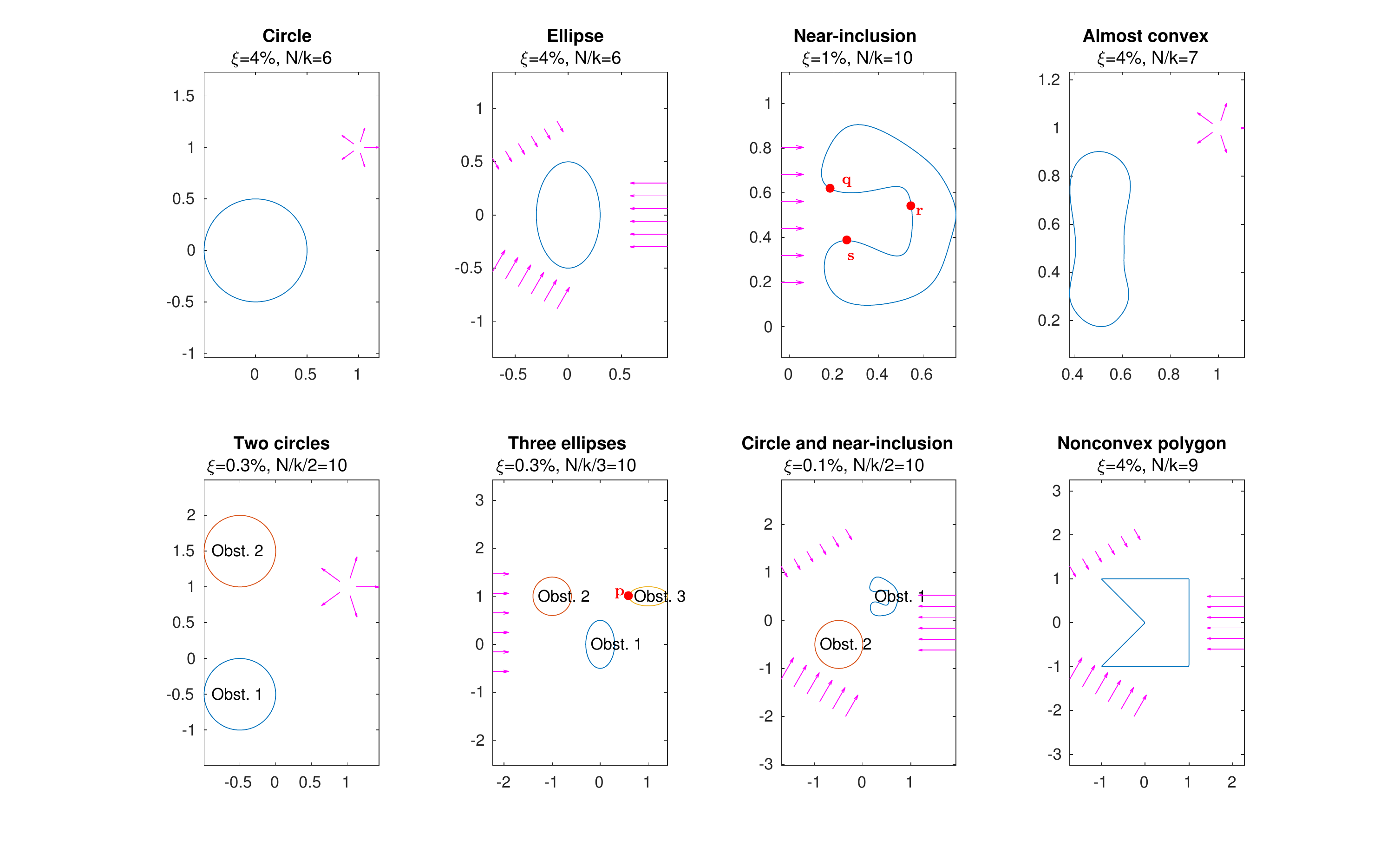} 
\caption{These are the obstacles used in the numerical experiments, with their respective names, threshold percentages $\xi$, relative number of points and incident waves. For the superposition of incident plane waves in different directions, the drawn sizes of the waves indicate their relative magnitude. Some points of relevance to the discussion are explicitly indicated in the figures and referenced in the text.}
\label{Fobsts}
\end{figure}

Throughout this section, we perform experiments in 2D using incident plane waves on Dirichlet obstacles shown and named in \cref{Fobsts}, where we also indicate the number of points per wavelength and the incident wave. These include a highly non-convex obstacle, multiple scattering configurations and a superposition of three incident plane waves with $2\pi/3$ angle differences. We will illustrate different aspects of our methods with these different scattering configurations.

We use piecewise linear basis functions in the reported experiments, since we explicitly aim to investigate asymptotic compression in low-order discretizations. The use of higher order basis functions leads to similar results in \cref{SaddRes}.

The scheme is implemented in Matlab, where we compute the integrals in \cref{EmatEntries} for computing the full matrix and \cref{Eslpot} for evaluation of the scattered field using the Sweldens quadrature scheme \cite{sweldens1994quadrature}. The BEM implementation that we started from is described in more detail in \cite{huybrechs2008gratings}. The integrals for $\tilde{A}_{m,n}$ and $R_{m,n}$ are replaced by their approximations \cref{eq:element_approximation,EroisSum}.

We use two error criteria in the examples. The first is the relative error in the solution vector $\tilde{c} = \tilde{A} \setminus b$ using the compressed matrix $\tilde{A}$ versus the exact numerical solution $c = A \setminus b$. The second criterion is the residual error of the integral equation, which for equation \eqref{Eie} corresponds to the accuracy with which the field satisfies the boundary condition. We have measured this error at $100$ points randomly distributed over the boundary of the scatterer and computed the relative 1-norm.

\subsection{Sparsity pattern and compression rate}

\begin{figure}[t]
\centering
\subfloat{\includegraphics[width=0.49\hsize]{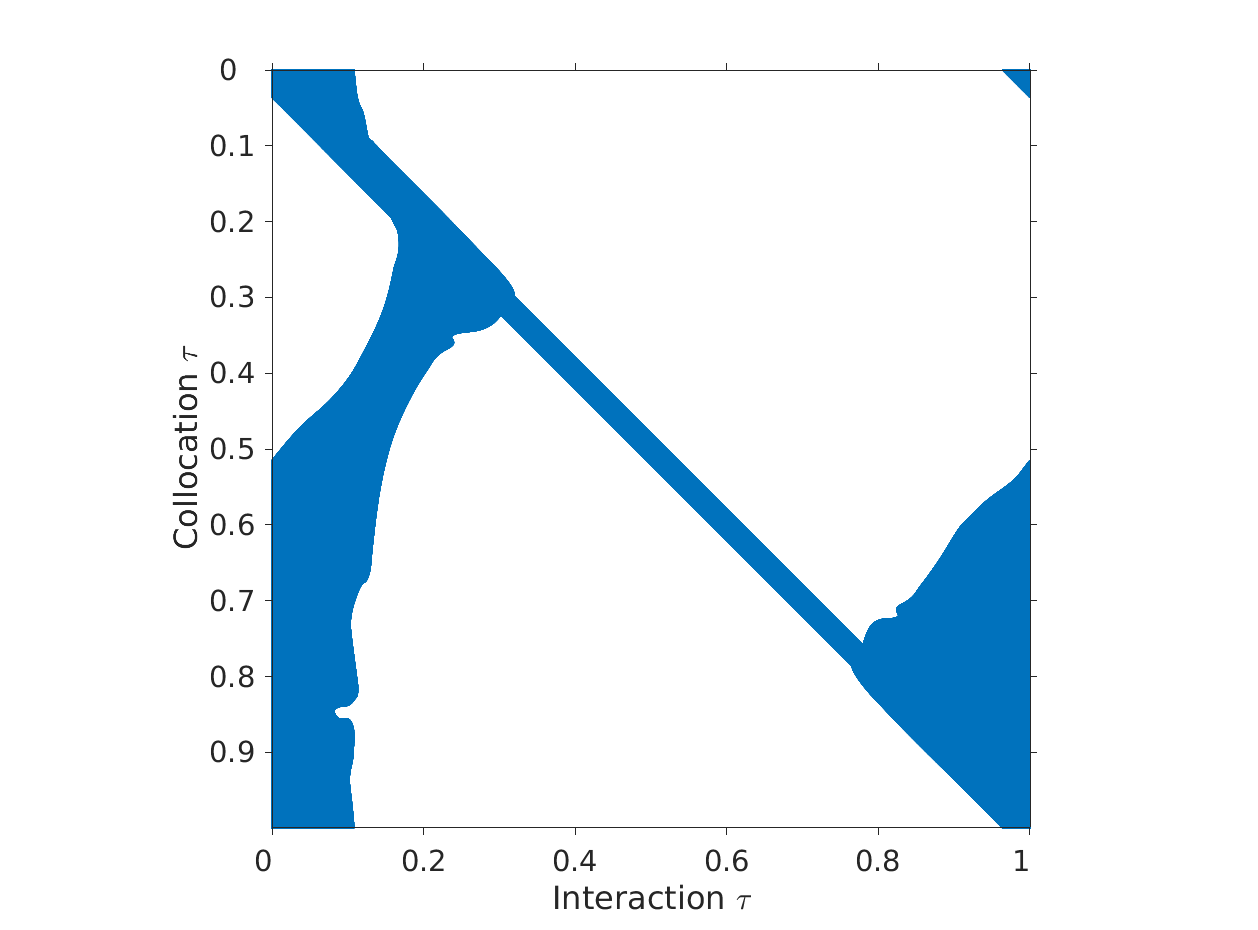} } 
\subfloat{\includegraphics[width=0.5\hsize]{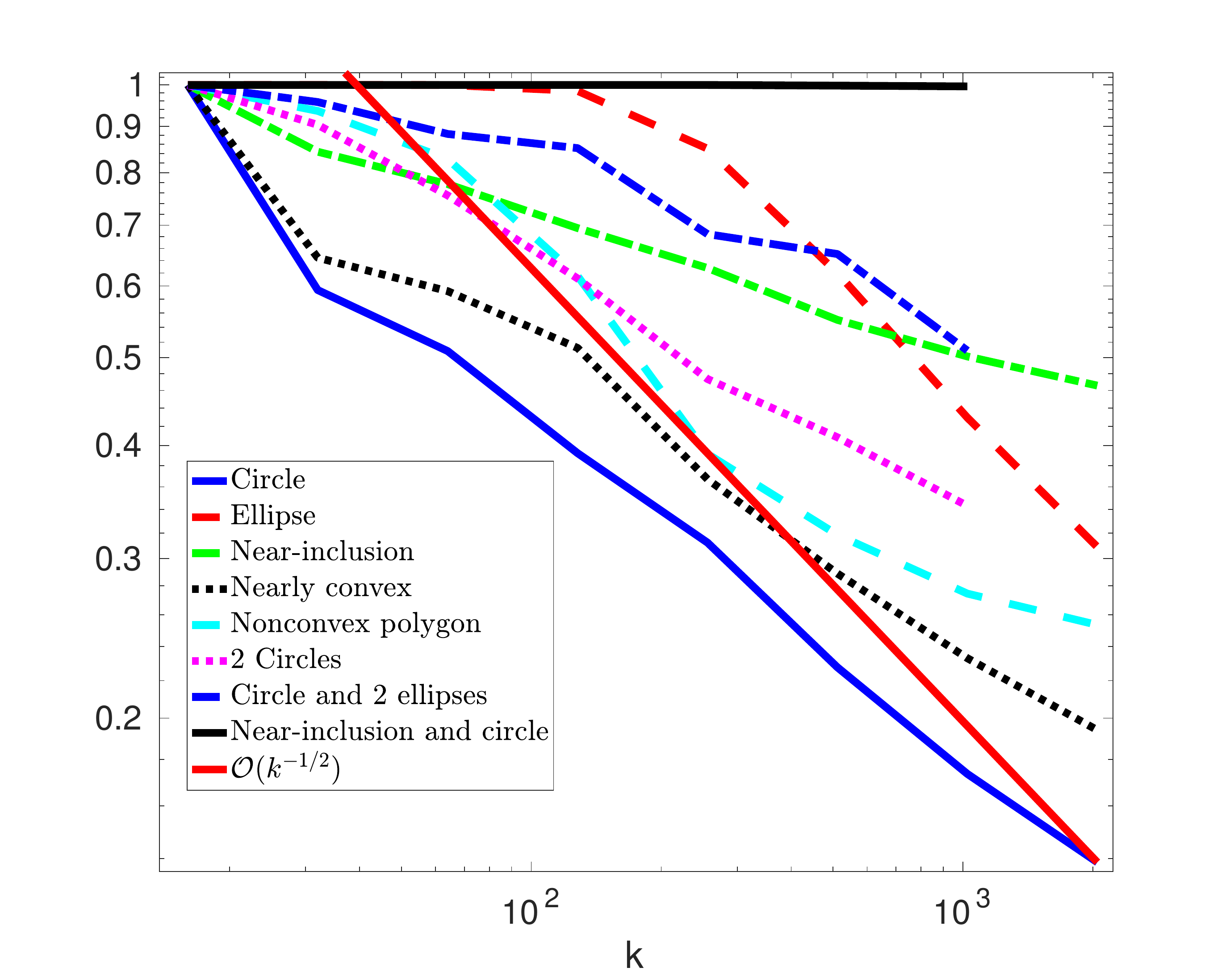} }  
\caption{Left: Structure of the compressed matrix $\tilde{A}$ using adaptive recompression on the almost convex obstacle at $k=2^{11}$. Right: fraction of nonzero elements in $\tilde{A}$ using adaptive recompression for increasing wavenumbers. The compression rate improves with increasing $k$ in all but one case and although still fairly modest, it is expected to keep improving for even higher wavenumbers, even for general obstacles with unknown high-frequency behaviour.}
\label{FrecomprPerc}
\end{figure}

We start by illustrating the sparsity of the discretization matrix $\tilde{A}$ and the compression rates for increasing $k$ for the obstacles shown in \cref{Fobsts}. We have used the adaptive recompression scheme outlined in \S\ref{SadapRecomprExplain}, starting from small values of the wavenumber and using threshold percentages and relative numbers of collocation points as indicated in \cref{Fobsts} along with each obstacle.

The sparsity pattern of $\tilde{A}$ is shown for the almost convex obstacle in the left panel of \cref{FrecomprPerc}. The diagonal clearly shows the contribution of the singularity of the Green's function on each row. There are two sidelobes, corresponding to collocation points in the shadow region having a stationary point in the illuminated region (see \S\ref{Sshadow}); however, these were determined automatically. The structure of this matrix resembles very much the structure of compressed matrices that appeared in literature for a strictly circular obstacle, for example in the asymptotic methods of \cite{sparseDiscr} and \cite{giladi}. Yet, there are two important main differences. On the one hand, the matrix in this paper is larger, since its dimensions scale linearly with increasing $k$. The methods of \cite{sparseDiscr} and \cite{giladi} used phase extraction to reduce this complexity, and as a result are limited to convex obstacles. On the other hand, the method of this paper does not require any a priori explicit knowledge of the nature of the solution. Therein lies a major gain. Though phase extraction is perhaps straightforward for the illuminated part of a convex obstacle, it is much less so on the shadow side. In the transitional region in between, the correct asymptotic behaviour becomes rather involved, even for convex obstacles \cite{unifPhaseExtr,MelroseTaylor}.

The right panel of \cref{FrecomprPerc} shows the improvement in compression rate as a function of the increasing wavenumber for all obstacles and it paints an encouraging 
picture. It is evident that the compression is complexity reducing for all but one configuration, including even the configurations with complicated incident waves, a near-inclusion and with multiple obstacles -- though at the shown range of the wavenumber, the percentage of non-zero elements remains fairly high. 

For the near-inclusion and circle multiple scattering configuration, there seems to be no compression. This is due to the combination of a complicated multiple scattering configuration and a superposition of three incident plane waves with $2\pi/3$ angle differences. The nonconvex polygon does not suffer as much from this, but the ellipse does at low wavenumbers. As there are $99.83\%$ nonzeros for the near-inclusion and circle at $k=512$ and $99.64\%$ at $k=1024$, visible compression is only expected to appear at much higher wavenumbers for this combination of a complicated geometry with a `difficult' boundary condition.

Yet, in spite of these issues, the example does illustrate that asymptotic information can be exploited even in the presence of highly non-trivial scattering patterns. Though compression may not result in large computational savings for challenging cases, later examples will show that the gain in condition number typically improves rapidly with increasing frequency and this may prove to be more important.

As expected, the results are best for convex and almost convex obstacles. These are the curves with decidedly better compression rates, corresponding to `easier' scattering problems: a circle, the almost convex problem and the non-convex polygon, and perhaps the ellipse and the two circles. For cases where asymptotics apply, asymptotic methods are likely superior, but our compression rate does seem to improve like $\MO(k^{-1/2})$. This corresponds closely to the expected scaling of $\MO(k^{-1/2})$ width of a window around a stationary point, as explained in \cref{s:integrals}. Indeed, because we refine the discretization with increasing wavenumber, there are $\MO(k)$ rows in the matrix. Each row has $\MO(k)$ elements, but each window covers only $\MO(k^{1/2})$ of them due to their decreasing support width. That leads to $\MO(k^{3/2})$ nonzero elements in $\tilde{A}$ in total, compared to $\MO(k^2)$ elements overall. Points in the transition region need larger $\MO(k^{-1/3})$ windows, but this region itself shrinks like $\MO(k^{-1/3})$ \cite{sparseDiscr}. Although this reasoning is a crude computation, it supports a compression rate of $\MO(k^{-1/2})$. This scaling seems to be matched nicely in the adaptive scheme. It is less clear what the compression rate in the other configurations is, as a function of the wavenumber. Yet, the downward trend is evident in all cases.

Thresholds were chosen to maintain roughly the same accuracy for increasing $k$, approximately equal to the discretization error. Precise results are given in \cref{ScompDiscErr} further on.

\subsection{The correlation matrix and adaptive localization of stationary points} 

Next, we illustrate the correlation matrix (introduced in \cref{Scorr}) and its relation to the structure of a scattering problem. As described earlier in the paper, first we solve the full scattering problem for some value of the wavenumber. The solution is used to estimate the relevance of parts of the boundary for each collocation point using a moving window in the oscillatory integral. In these experiments, we have used the value $T=0.02$ for the parameter that determines the width of the moving window.

\begin{figure}[h]
\centering
\subfloat{\includegraphics[width=0.5\textwidth]{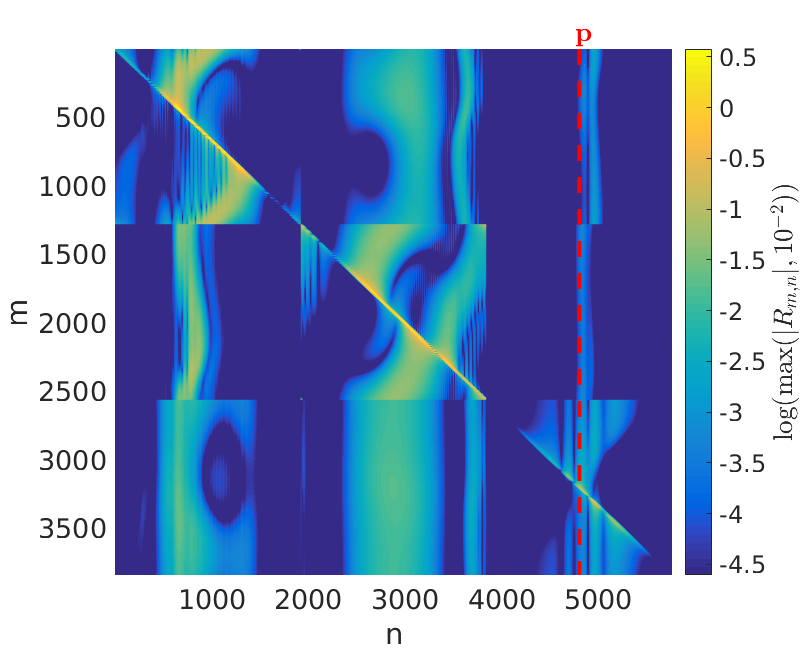} }
\subfloat{\includegraphics[width=0.5\textwidth]{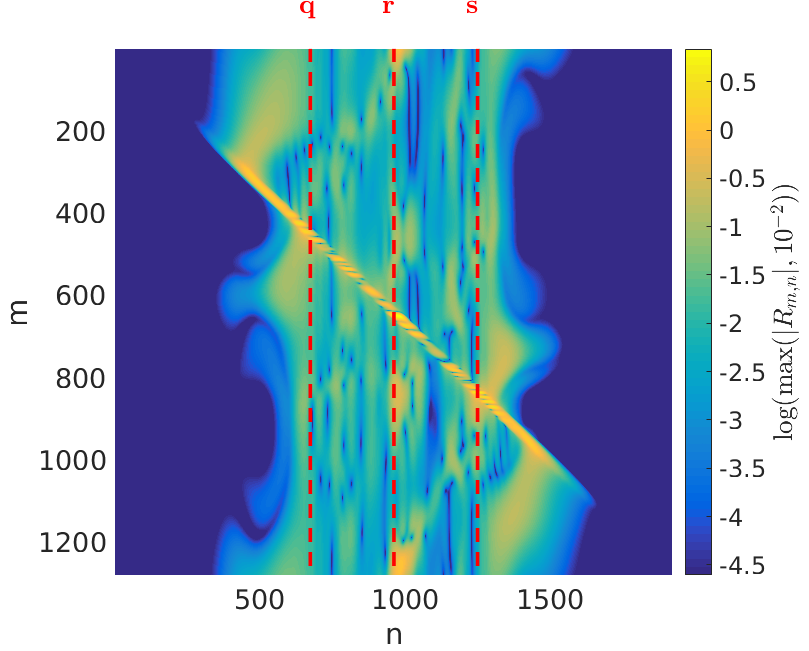} }
\caption{Correlation matrix $R_{m,n}$ at $k=2^7$ for the multiple scattering configuration with three ellipses (left), and for the scattering configuration with the near-inclusion obstacle (right). Each row corresponds to a collocation point $\tau_m$ on one of the obstacles. Each column corresponds to a center of a moving window $\sigma_q$. 
The value shown is the value of the localized oscillatory integral, i.e. of the oscillatory integral operator applied to a windowed part of the exact solution through the approximation \cref{EroisSum}.}
\label{Froisspy}
\end{figure}

An informative example of a correlation matrix is shown in the left panel of \cref{Froisspy}, for the multiple scattering configuration with three ellipses (the ellipses are shown and numbered $1-2-3$ in \cref{Fobsts}). Here, rows correspond to collocation points and columns to the centers of the moving windows. The correlation matrix has a $3 \times 3$ blocked structure, with each block corresponding to the coupling between two of the three ellipses. Asymptotic analysis of multiple scattering configurations shows that ultimately, rays bounce around in periodic orbits that mimize the distance between different obstacles~\cite{2DEcevit,3DEcevit}. These periodic orbits appear as the vertical patterns in the coupling submatrices in \cref{Froisspy}.

For example, consider the two blocks in the upper right of the computed correlation matrix, on the first and second rows. These two blocks correspond to collocation points on ellipse $1$ or ellipse $2$, with an observation point varying on the third obstacle. The center point $\mathbf{p}$ that is indicated in the left part of \cref{Froisspy} (and in \cref{Fobsts}) corresponds to the leftmost point on the third ellipse. The correlations are high when concentrated near $\mathbf{p}$, and this is physically the only region where rays can reflect off the third ellipse onto the other two along the shortest path, for a plane wave incidence from the left. Hence, that is where any of the stationary points must lie. The compressed matrix $\tilde{A}$ that follows from this $R_{m,n}$ can be rather sparse away from $\mathbf{p}$ in these two subblocks.

The Green singularity is clearly visible in the correlation matrix along the diagonal. Note that the three diagonal blocks are the self-interaction of the ellipses. These display a pattern that relates to that shown in \cref{FrecomprPerc} before. The first and second diagonal block exhibit large correlations in much wider anti-diagonal side-lobes above and below the main diagonal. These correspond to a stationary point in the illuminated region, for a collocation point in the shadow region of the obstacles, much like for the single scattering problem shown in \cref{FrecomprPerc}. This is less pronounced on the third obstacle, because it lies in the shadow of the first one.

The results are similar for the non-inclusion single scattering obstacles. More specifically, a collocation point in the shadow region has a stationary point in the illuminated region, in agreement with \cref{Sshadow}. The computed correlations as in \cref{Froisspy} seem quite good at following the stationary points that one can deduce from \cref{Sshadow} or reflection arguments. 

For the near-inclusion obstacle in the right part of \cref{Froisspy}, $R_{m,n}$ is high for all $n \in [600,1300]$ for most $m$. So in $\tilde{A}$, almost all of the illuminated region is included in the windows, due to the trapping of rays inside the cavity. As a result, there is very little compression. Of course, one could not hope for much better in view of the complicated high-frequency behaviour of this obstacle, which figures like these could help to understand. Additionally, there is a very high correlation near the points $\mathbf{q}$, $\mathbf{r}$ and $\mathbf{s}$ which are also indicated in \cref{Fobsts}. An important side-effect of the complicated high-frequency behaviour is that some of the correlations near the diagonal, covering the Green function's singularity, were rather small. In the implementation we always include the Green singularity into the support of the window function, regardless of the value of the correlations.

\begin{figure}[h]
\centering
\subfloat{\includegraphics[width=0.5\hsize]{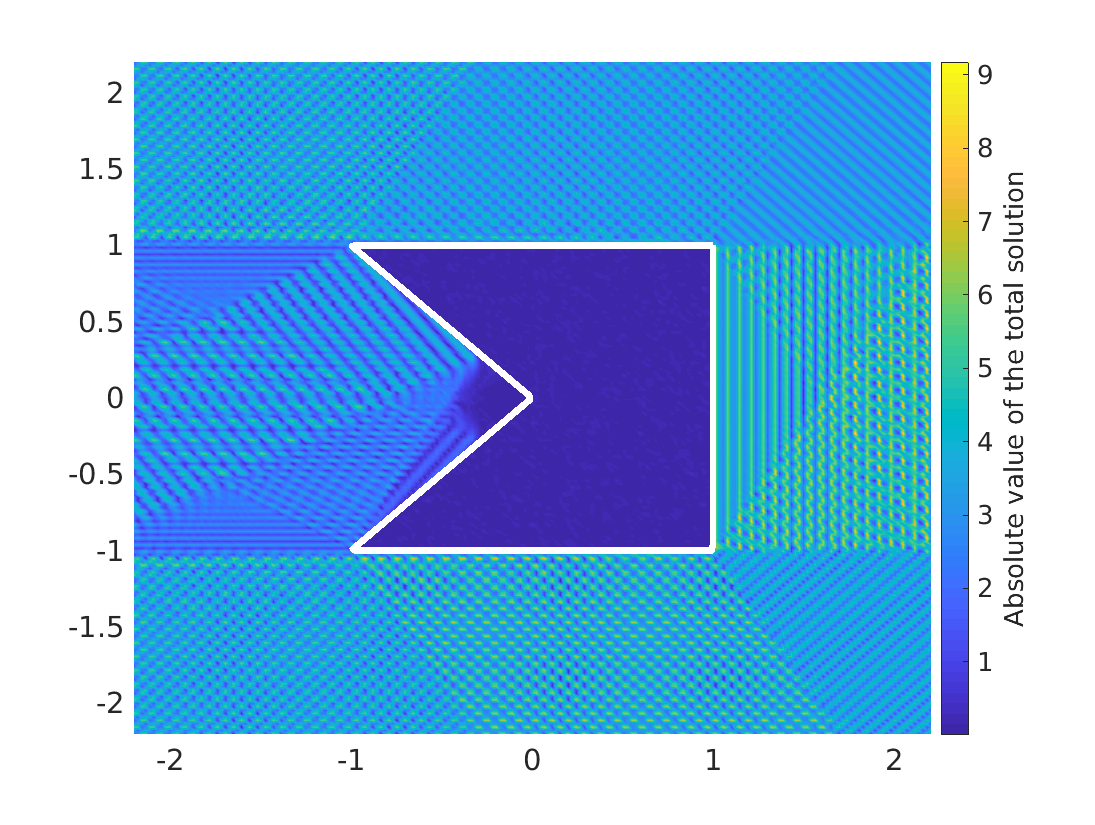} } 
\subfloat{\includegraphics[width=0.5\hsize]{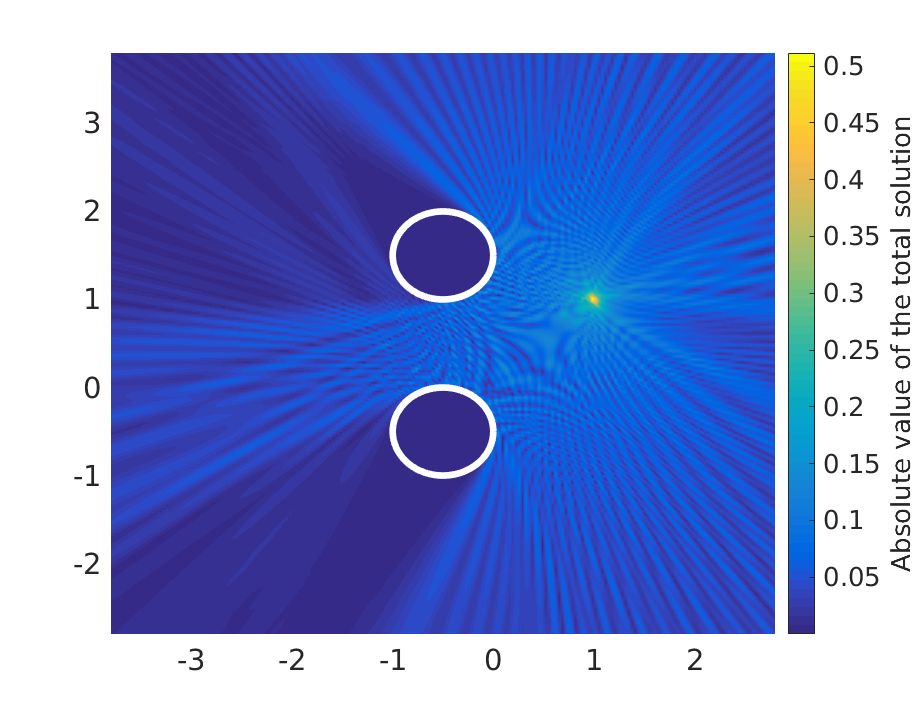} }
\caption{Total field for the nonconvex polygon at $k=2^{11}$ (left) and for the two circles at $k=2^7$.}
\label{Ffield}
\end{figure}

One may understand these oscillations of the correlations near the diagonal by viewing the total field obtained after solving the system. The simple example at the left of \cref{Ffield} shows an interference pattern within the nonconvex part of the obstacle. Although the near-inclusion obstacle exhibits significantly more complicated high-frequency scattering behaviour, the nonconvex polygon also has some small correlations near the diagonal, so we expect these to originate in interference. Similarly, there is interference between the two circles shown in the right of \cref{Ffield}. This is the case for the lower right part of the upper circle and the upper right part of the lower circle, as they can recieve rays from the point source at $(1,1)$ as well as from the other circle.

\subsection{Compression error and discretization error} \label{ScompDiscErr}

\begin{table}[h]
\centering
\caption{Accuracy as a function of the wavenumber for the scattering configuration with two circular obstacles. The windows were computed using correlations at $k=128$ and reused for larger wavenumbers. The errors indicate how well the scattered field matches the boundary condition, for the exact solution $c = A^{-1} b$ and the approximate solution $\tilde{c} = \tilde{A}^{-1} b$ computed after asymptotic compression.}
\label{TcorrMOerrBC}
\begin{tabular}{l|cccc}
\hline\noalign{\smallskip}
Residue int. eq. & $k=128$ & $k=256$ & $k=512$ & $k=1024$ \\ 
\noalign{\smallskip}\hline\noalign{\smallskip}
$c$ & 6.78e-4 & 4.82e-4 & 2.55e-4 & 2.34e-4 \\ 
$\tilde{c}$ & 7.10e-4 & 4.86e-4 & 2.56e-4 & 2.34e-4 \\ 
\noalign{\smallskip}\hline
\end{tabular}
\end{table}

The overall goal is to match the error due to asymptotic compression to the discretization error, as one can not expect to improve on that. The correlation thresholds were chosen accordingly. The accuracy results for the multiple scattering configuration with two disks are summarized in \cref{TcorrMOerrBC}. Here, the discretization error is around $0.02\%$ using the number of points per wavelength indicated in \cref{Fobsts} in our standard BEM implementation. Using a threshold percentage of $\xi = 0.3\%$, we can compute the scattered field using the compressed matrix $\tilde{A}$ without additional loss of accuracy. In this example, we have used a fixed window computed at $k=128$ and reused that window for larger values of $k$.

\begin{table}[h]
\centering
\caption{Results similar to those in \cref{TcorrMOerrBC}, but for the near-inclusion obstacle and using adaptive recompression. An advantage of recompression is that we can start at a smaller value of the wavenumber.}
\label{TrecomprErrBC}
\begin{tabular}{l|cccccccc}
\hline\noalign{\smallskip}
\begin{minipage}{0.08\textwidth} Residue\\ int. eq. \end{minipage} & $k=16$ & $k=32$ & $k=64$ & $k=128$ & $k=256$ & $k=512$ & $k=1024$ & $k=2048$ \\ 
\noalign{\smallskip}\hline\noalign{\smallskip}
$c$ & 1.81e-2 & 1.65e-2 & 1.32e-2 & 1.39e-2 & 1.16e-2 & 1.62e-2 & 1.59e-2 & 1.87e-2 \\ 
$\tilde{c}$ & 1.81e-2 & 1.75e-2 & 1.42e-2 & 1.45e-2 & 1.19e-2 & 1.65e-2 & 1.61e-2 & 1.90e-2\\ 
\noalign{\smallskip}\hline
\end{tabular}
\end{table}

The results when using adaptive recompression at higher frequencies are entirely similar, since the recompressed correlation matrix $R_{m,n}$ at a wavenumber $k_j$ behaves in the same way as when we had computed full correlations at $k_{j}$: we just recursively leave out negligible contributions. One example is shown in \cref{TrecomprErrBC} for the near-inclusion obstacle. This example shows that the errors of the method do not deteriorate when using adaptive recompression, even for an obstacle with complicated high-frequency behaviour.

\begin{figure}[t]
\centering
\includegraphics[width=0.7\textwidth]{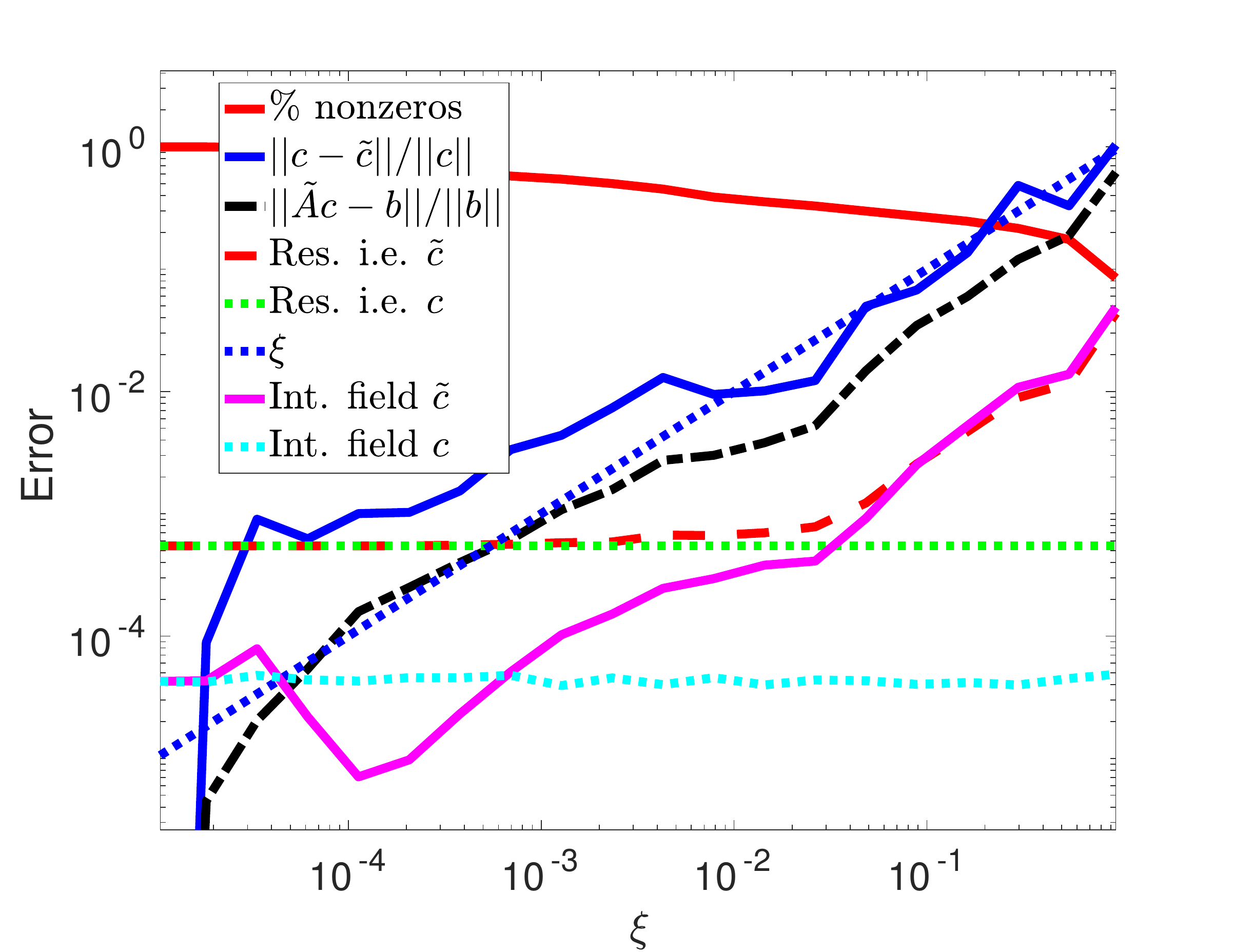} 
\caption{Fraction of nonzeros in $\tilde{A}$ and errors as a function of the threshold percentage $\xi$ at wavenumber $k=2^8$ for the almost convex obstacle. The errors improve with increasing $\xi$, but the sparsity decreases.}
\label{Fthr}
\end{figure}

\Cref{Fthr} shows the percentage of nonzero elements as a function of the threshold percentage $\xi$, for the almost convex obstacle. For low $\xi$, $\tilde{c}$ provides a good approximation of the solution vector, as our compressed matrix $\tilde{A}$ has a low forward compression error then. The errors on the scattered field also decrease with $\xi$ until they reach the order of the error when using the full matrix $A$. This is shown in two ways: the error in satisfying the boundary conditions (the residual error of the integral equation), and the total interior field. The latter should be zero for this particular Dirichlet problem and it saturates at about $0.004\%$ with and without compression, averaged over $100$ points inside the obstacle chosen randomly for each $\xi$.

\subsection{Computational efficiency and condition numbers} 

\begin{figure}[h]
\centering
\includegraphics[width=0.7\textwidth]{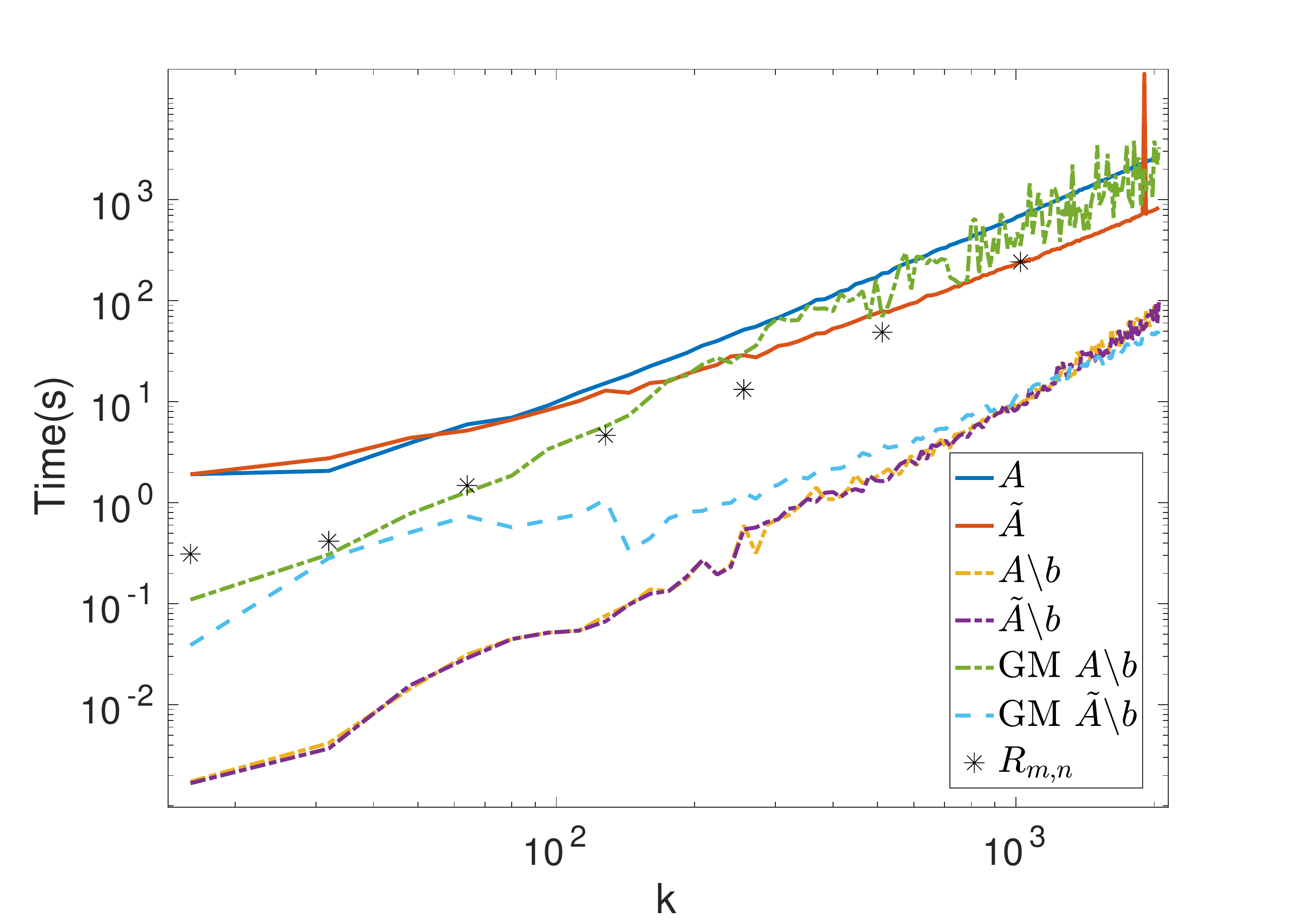}
\caption{Time (s) at each wavenumber for the different computations needed for the nonconvex polygon, where we can choose between gmres and the \textsc{matlab} backslash. We use linear basis functions and recompression in $\tilde{A}$. The wavenumber $k$ increases linearly in steps of $16$, while our correlations are efficiently recomputed for $k$ at integer powers of two.}
\label{TmoRecTime}
\end{figure}

Now, we summarize performance aspects of the adaptive asymptotic recompression scheme that was introduced in \cref{SadapRecomprExplain}. First, we measure the efficiency of the frequency sweep with adaptive compression, compared to one where the full discretization matrix is constructed. We have used Matlab2017a on a $64$-bit machine with $66$ GB memory and $32$ Intel(R) Xeon(R) CPU E5-2650 v2 CPU's at $2.6$ Ghz: the results are shown in \cref{TmoRecTime} for the nonconvex polygon. 
Although no other tasks were burdening this machine, there is an outlier at $k=1904$. At each other wavenumber, it is faster to construct $\tilde{A}$ than the full matrix $A$, though the factor is not very large in this example. However, a large gain is observed when one would solve the resulting linear system using an iterative solver. Moreover, the gains accumulate in a frequency sweep.

For this example, the total cost summed over the $128$ linearly spaced (as in an FFT) wavenumbers is $1.20e5$ seconds. This consists in constructing a full $9k \times 9k$ matrix and solving the associated linear system using a direct solver at wavenumber $k_1=16$, then $k_2 =32$ and so on until $k_{128} = 2048$. In contrast, the total cost of computing all correlations, constructing all $128$ compressed matrices and solving the resulting linear systems, is 4.04e4 seconds excluding the outlier, or one third of the cost of the full method. This gap is only expected to increase for larger wavenumbers as the compression rates will improve. Especially when one would use iterative solvers, solving the linear system will become the dominant cost, which our method decreases through a sparse matrix and low number of iterations.

\begin{figure}[h]
\centering
\subfloat{\includegraphics[width=0.5\hsize]{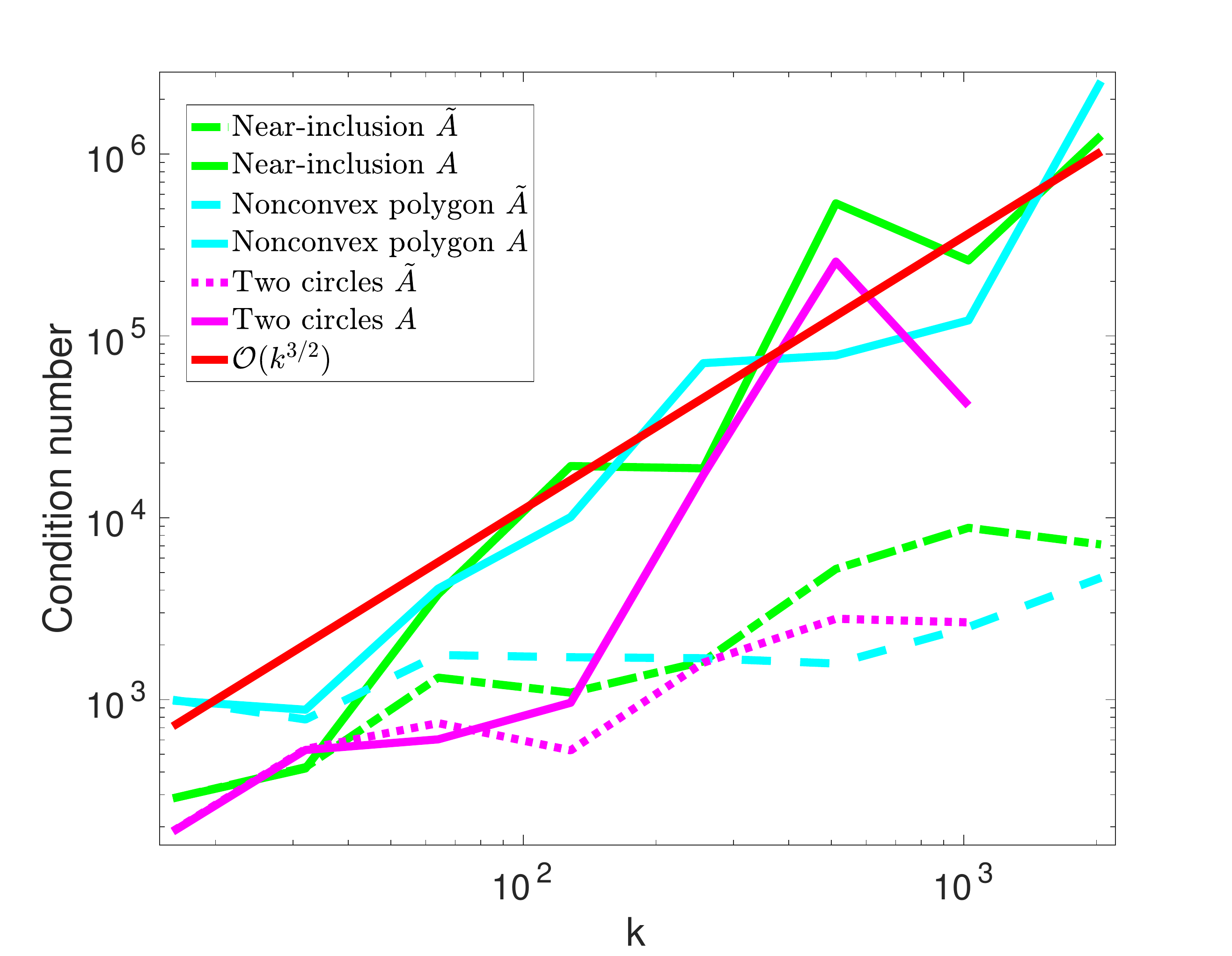} } 
\subfloat{\includegraphics[width=0.5\textwidth]{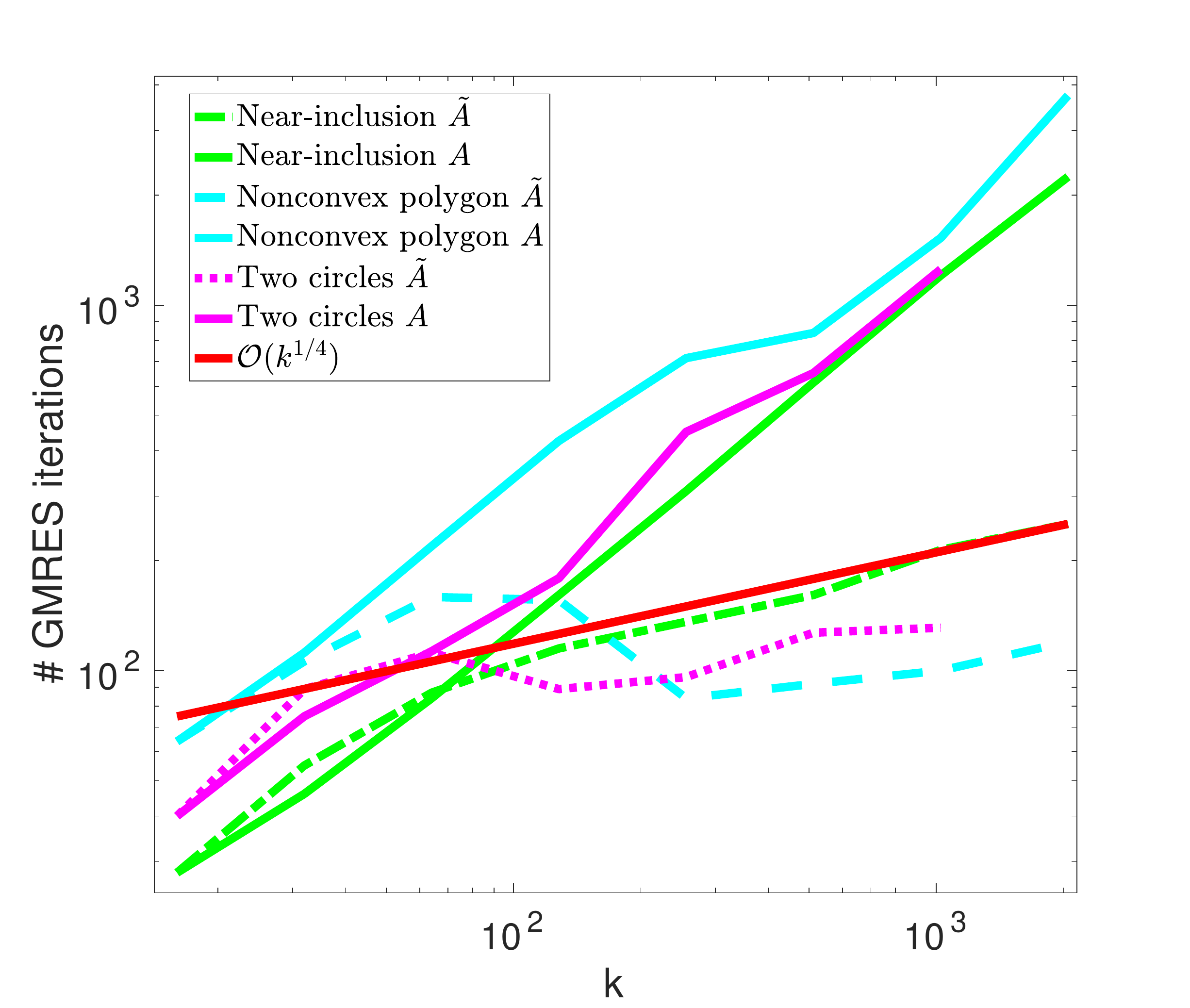} } 
\caption{Estimated condition number and number of GMRES iterations for some obstacles when using adaptive recompression and without compression for wavenumbers at integer powers of two.}
\label{FcondGmRecompr}
\end{figure}

Finally, \cref{FcondGmRecompr} shows that the compression leads to an improvement in condition number, which then lowers the number of iterations of an iterative solver. In general, however, the obstacle and/or the incident wave may be too intricate to allow any significant compression or such speedup of the iterative solution at the wavenumbers considered. We have used a non-restarted GMRES in Matlab with a tolerance of $10^{-5}$ for this experiment. The same qualitative behaviour is observed for all obstacles, using piecewise linear basis functions. The condition number improves by a factor of $1.1$ to $3$ compared to using fixed windows from correlations computed at $k=2^7$, although the number of iterations varies less predictably about that range. % Maybe remove last sentence.

\subsection{Visibility criterion} 

\begin{figure}[t]
\centering
\subfloat{\includegraphics[width=0.49\hsize]{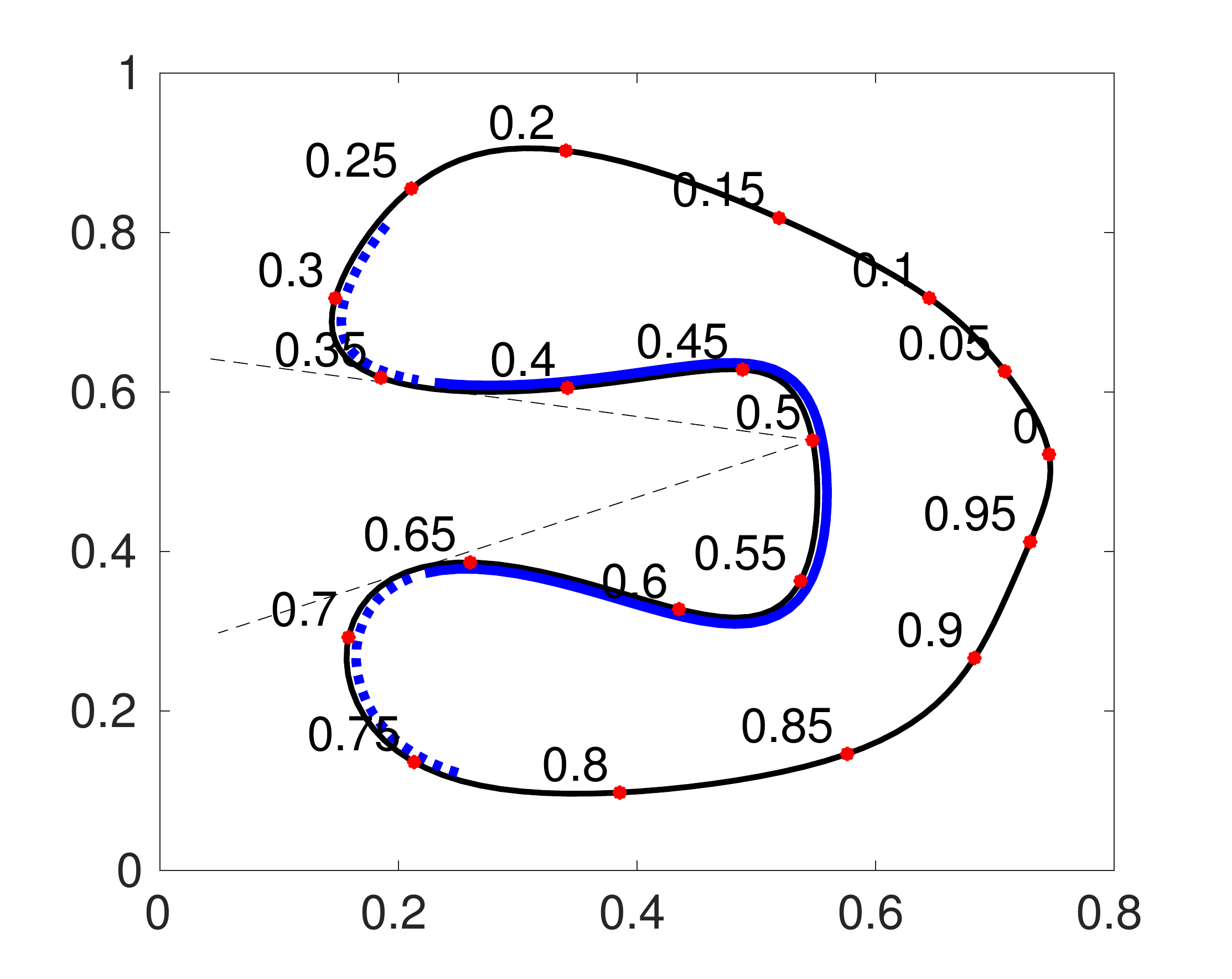} } 
\subfloat{\includegraphics[width=0.5\hsize]{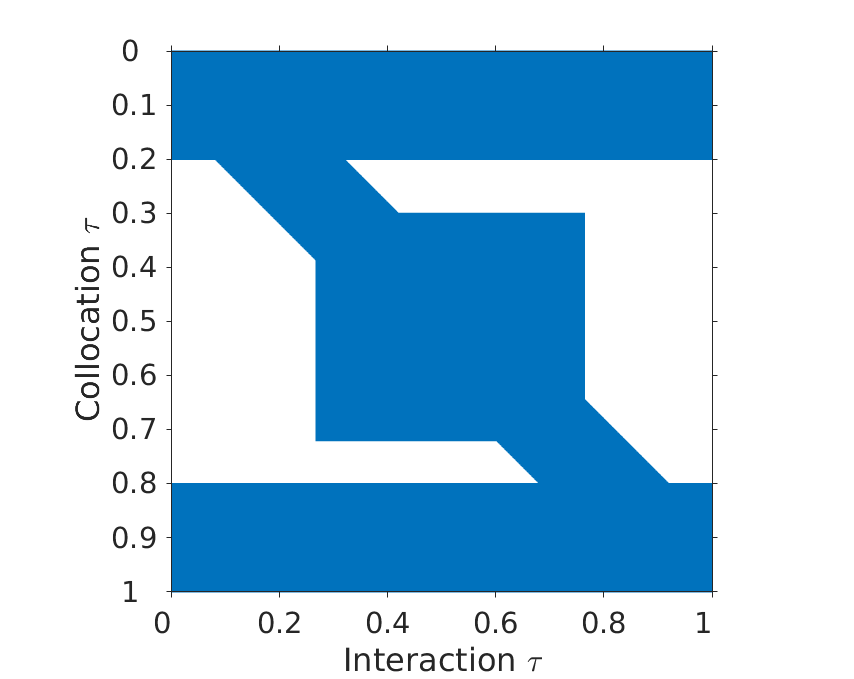} }
\caption{Left: Obstacle, parametrisation and region which is visible from $\tau=1/2$ for the near-inclusion obstacle. Right: Structure of the compressed matrix using the visibility criterion.}
\label{FvisCone}
\end{figure}

Next, we illustrate the possibilities and limitations of compression using just the visibility criterion outlined in \cref{SvisExplain}. We focus on scatterers with complicated asymptotic structure of the rays, since one can not gain much from asymptotics in such cases. We use the near-inclusion domain as a first example. The domain and its parameterization are shown in \cref{FvisCone}. The corresponding sparsity pattern of $\tilde{A}$ is shown in the right of \cref{FvisCone}, where we have used the visibility criterion to perform compression only for collocation points in the illuminated region.

We have highlighted a collocation point corresponding to $\tau = 0.5$, and \cref{FvisCone} shows the part of the obstacle that is directly visible from this point. It includes much of the cavity, but nothing outside of it. We add a dashed blue region where the window can smoothly decay to zero. This window translates into the large dense square in the centre of the sparsity pattern of $\tilde{A}$: all points inside the cavity are asymptotically linked directly to each other. The dense block captures the complicated multiple scattering going on inside the cavity. Yet, this whole block is asymptotically independent of the solution outside the cavity, and this yields the zero elements outside the range $[0.365-T, 0.665+T]$ for collocation $\tau \in [0.3,0.72]$, where the length of the decaying part $T=0.1$. 

Collocation points with $\tau$ approximately in $[0.2,0.3]$ or $[0.72,0.8]$ correspond to a locally convex part of the obstacle in the illuminated region. These points would only be influenced by themselves, as they cannot `see' other points. This results in a small band $[\tau-0.02-T, \tau+0.02+T]$ around the diagonal in the discretization, above and below the central square at the right of \cref{FvisCone}. As mentioned, we have not attempted to apply compression based on visibility in the shadow region ($\tau \notin [0.2,0.8]$). 

We compute $\tilde{A}$ at $k=2^9$ using this structure, which results in a relative error on $c$ of $12.1\%$ and a worsening of the residue of the integral equation on 100 random points from $1.19\%$ to $5.84\%$. This gives a qualitatively accurate solution, but the errors are quite large. The errors do respectively improve to $3.85\%$ and $2.38\%$ by not applying compression in the locally convex parts of the obstacle in the illuminated region. Or one can increase the frequency, as it seems that the inclusion region still significantly influences those parts at $k=2^9$ due to creeping rays.

One may expect that the visibility criterium might be easier to validate for multiple scattering of convex obstacles, or obstacles with a finite number of reflections. Here, one has to be careful in what it means to be `visible'. We test the multiple scattering configuration with two circles shown in \cref{Fobsts} at $k=256$ where we change the point source to a plane wave coming from the left to obtain the more pronounced behaviour described below. First, we compress parts of the rows that correspond to self-interactions in the region that is directly illuminated by the source. More precisely, if the $x$-coordinate is less than $-0.5-T$ at $\kappa(t)$, then $\lambda = t-2T, l = t-T, r = t+T$ and $\rho = t+2T$ for $T=0.15$ and there is no compression elsewhere. However, this worsens the relative residue of the integral equation from $1.0835\%$ 
to $5.41\%$ and gives a relative error on $c$ of $21.28\%$.

\begin{figure}[h] 
\centering
\subfloat{\includegraphics[width=0.65\textwidth]{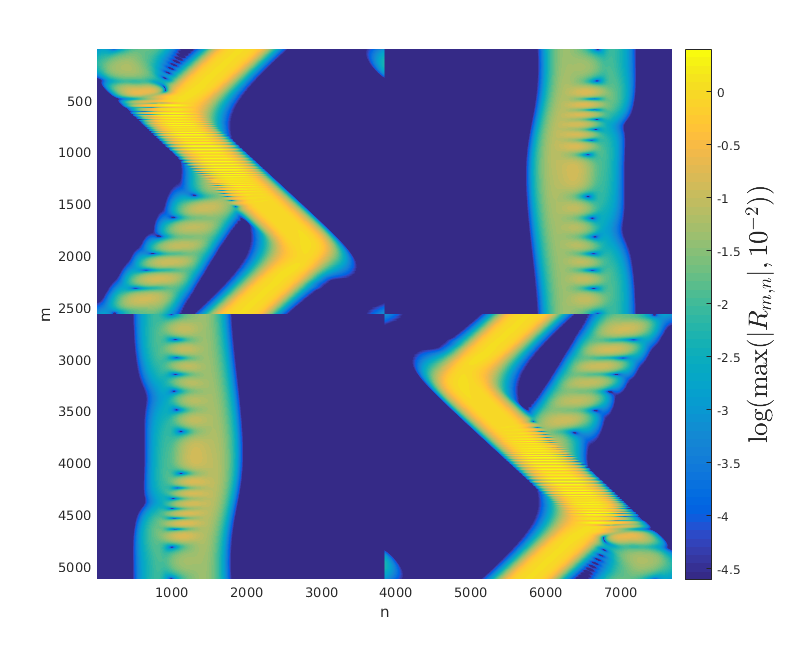}} 
\subfloat{\raisebox{1.4cm}{\includegraphics[width=0.4\textwidth]{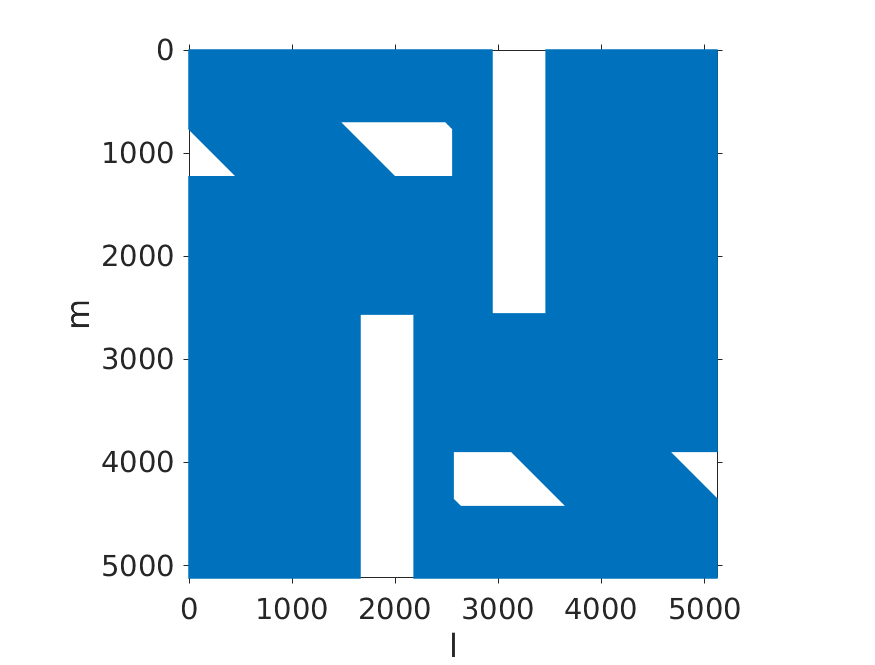}}}
\caption{Correlations for the two circles with an incident plane wave from the left at $k=256$ (left) and sparsity pattern using visibility (right).}
\label{FcorrTwoCirc}
\end{figure}

These high errors are explained by the correlation matrix shown in \cref{FcorrTwoCirc}. In the middle rows of the self-interaction blocks, one not only notices the Green singularity on the diagonal, but also a stationary point elsewhere. For example, in the upper left block of \cref{FcorrTwoCirc}, there is a side lobe along the antidiagonal direction that ranges approximately from the row labelled $1500$ to the one labelled $2500$. These rows correspond to collocation points in the lower half of the lower obstacle. These points are in the shadow with respect to the upper obstacle, and therefore they have a stationary point in the illuminated part, i.e. the upper half of the lower circle. The concepts of illumination and shadow are relative to the other obstacle here, not to the incoming plane wave. 

Applying compression solely for collocation points on parts of the obstacle that are illuminated both by the incoming wave and by the other obstacle is successful: the residue of the integral equation and the error on $c$ become $1.0832\%$ and $0.0816\%$ respectively. This indicates that each obstacle in a multiple scattering configuration should be viewed as a source on top of the incoming wave when applying the visibility criterion.

Yet, much more compression is possible: it is quite clear from the correlation matrix in \cref{FcorrTwoCirc} that there is large potential for compression in the coupling blocks. Indeed, any point in the upper half of the upper obstacle can not `see' the whole lower obstacle. Similarly, no point in the lower half of the lower obstacle can be a stationary point for a collocation point anywhere on the whole upper obstacle. In the upper right and lower left blocks of the correlation matrix, this yields large blue regions amenable to compression. The additional compression leads to the $85\%$ nonzeros in the right of \cref{FcorrTwoCirc}, a residue of the integral equation of $1.0839\%$ and an error on $c$ of $0.0817\%$. The compression rate could improve with increasing frequency as the width of the decaying part $T$ can decrease.

The visibility approach could be made automatic and more robust, and whether two points can `see' each other is independent of the incident wave. However, one can only gain a percentage of nonzeros that is independent of $k$ and which deteriorates with more complex geometries or more obstacles, as each of them acts as a source giving a shadow. The shadow region is also dependent on the incident wave and we applied no compression there. The visibility criterion would benefit from a better understanding of the stationary points of collocation points that are in the shadow (recall \cref{Sshadow}).

\subsection{Additional robustness results} \label{SaddRes}

To begin, compressing the BEM matrix also helps against ill-conditioning of the integral formulation of the first kind near resonances. As a simple example, the resonant frequencies for our circle of radius $1/2$ are $k=2j_{n,m}$ (the $m$th zero of the $n$th unmodified Bessel function of the first kind). The discretization matrix of this problem in our code has condition number $1020$ for $n=50$ and $m=1$. The condition number is large but not infinite, due to the discretization error. Still, after adaptive compression with the point source from \cref{Fobsts}, the condition number is reduced to $36.9$ with $34\%$ nonzeros, while for an incident plane wave from the left, the condition number is reduced to $37.3$. The eigenfunction of the problem is inherent to the domain, not to the boundary condition, and it seems entirely eliminated by the compression.

Secondly, we illustrate a simplification of the compression scheme based on hard truncation: we simply discard all elements of $A$ where they would be zero in $\tilde{A}$ in the adaptive scheme, but we make no modifications to the other entries of $A$. Indeed, this corresponds closely to using a block window function, rather than the smooth function~\eqref{Ewin2D}. As mentioned in \cref{s:integrals}, this should lead to spurious $\MO(k^{-1})$ endpoint contributions for a general oscillatory integral. We consider the near-inclusion obstacle at $k=128$. Without compression, there is a residue of $1.09\%$ on the integral equation; with smooth windows, this becomes $1.26\%$ and $2.80\%$ error on $c$; with the non-smooth compression windows, only $3.33\%$ and $12.1\%$ respectively. The errors clearly grow, but not excessively so. The reason is, at least in this case, that the contributions of the Green's singularity and stationary points dominate. The spurious contributions caused by the discontinuous compression scheme are more than an order of magnitude smaller, such that one still obtains the correct qualitative behaviour of the solution.

Finally, the previous experiment is repeated with piecewise cubic basis functions $\varphi_j(\tau)$. The condition number improves from 1.28e3 to 295, while this was 1.27e3 to 116 for linear basis functions. The residue on the integral equation marginally changes from $1.16\%$ to $1.33\%$, for an error on $c$ of $4.66\%$. Although we mainly considered piecewise linear basis functions in this article, adaptive asymptotic (re)compression thus applies to more general basis functions.

\section{Concluding remarks} \label{Sconc}

Our results show that exploiting asymptotic behaviour in scattering problems is possible without a priori knowledge of the asymptotic behaviour of the solution. Compared to true asymptotic methods, our approach is costly as all oscillations are still resolved with a fine discretization. Yet, the compression scheme applies to general objects, arbitrary incident waves and appears to be robust and automatic in several ways for problems with highly complicated high-frequency asymptotics. In particular, we showed that asymptotic compression can be applied simply by rescaling standard BEM matrix elements, even in the presence of near-cavities.

There are several opportunities for further research. First, though we have avoided it in this paper, the present method may be combined with asymptotic methods. Any available asymptotic information can be used to reduce the cost of computing the solution. One can combine approximate asymptotic knowledge with the robustness of the present approach, for example if a small number of geometrically complicated parts of the boundary prevent a full asymptotic scheme. Figures like \cref{Froisspy} could help obtaining this asymptotic knowledge for intricate high-frequency behaviour.

Alternatively, the method can be combined with non-asymptotic methods. The Fast Multipole Method yields a fast matrix vector product for discretizations of scattering problems, with optimal time and memory complexity. Since our method essentially amounts to modifying the Green's function, the two methods may be compatible. 
Both may also be combined with local quadrature correction or quadrature by expansion \cite{qbx}, which are expected to improve the accuracy of the integrals we compute.

Finally, the principles underlying the asymptotic compression scheme are valid in other integral equation formulations, higher dimensions and other physical problems as well. Preliminary investigations using the 3D Galerkin scheme in the C++ software package BEMPP \cite{bempp} support this future research direction.

%\begin{acknowledgements}
\section*{Acknowledgements}
The authors would like to thank Samuel Groth, Stephen Langdon, Niels Billen, Philip Dutr\'{e}, Karl Meerbergen, Laurent Jacques and Dave Hewett for interesting and helpful discussions on this paper and related topics. The authors were supported by FWO Flanders [projects G.0617.10, G.0641.11 and G.A004.14].
%\end{acknowledgements}

\bibliographystyle{abbrv}
\bibliography{asyComprV19arXiv.bbl}

\end{document}